\documentclass[a4paper]{article}

\usepackage[pdftex,colorlinks,bookmarksnumbered,bookmarksopen,citecolor=red,urlcolor=red]{hyperref}
\usepackage[dvipsnames]{xcolor}
\usepackage{amsmath, amsthm, amsfonts, bbold, float, tikz, ulem, mathtools, tabularx, mathrsfs, booktabs}
\usetikzlibrary{shapes, decorations.pathreplacing, patterns}
\usepackage{multirow}
\usepackage{subfigure}
\usepackage{comment}

\newcommand{\bm}[1]{\text{\boldmath $#1$\unboldmath}}
\newcommand{\bx}{\bm{x}}

\newcommand{\bn}{\bm{n}}

\newcommand{\balpha}{\bm{\alpha}}

\newcommand{\bmu}{\bm{\mu}}

\newcommand{\bLambda}{\bm{\Lambda}}

\newcommand{\pgd}{\texttt{PGD}}

\newcommand{\nAIP}{N_{\texttt{AIP}}}
\newcommand{\nDP}{N_{\texttt{DP}}}
\newcommand{\nIP}{N_{\texttt{IP}}}
\newcommand{\nGMRES}{N_{\texttt{GMRES}}}
\newcommand{\dIP}{d_{\texttt{IP}}}
\newcommand{\Toff}{T_{\texttt{off}}}
\newcommand{\Ton}{T_{\texttt{on}}}

\newcommand{\hOmega}{\hat{\Omega}}
\newcommand{\Map}{\bm{\mathcal{M}}_{\mu}}
\newcommand{\hx}{\hat{x}}
\newcommand{\hy}{\hat{y}}
\newcommand{\hmu}{\hat{\mu}}

\newcommand{\mat}[1]{\mathbf{#1}}

\newtheorem{tm}{Theorem}
\theoremstyle{remark}
\newtheorem{rem}[tm]{Remark}

\pgfdeclarepatternformonly{mydots}
  {\pgfpoint{0pt}{0pt}}          
  {\pgfpoint{2pt}{2pt}}          
  {\pgfpoint{2pt}{2pt}}          
  {
    \pgfpathcircle{\pgfpoint{1pt}{1pt}}{0.4pt}
    \pgfusepath{fill}
  }

\usepackage{xifthen}
\newcommand{\pateraBlock}[6]{

    {\ifthenelse{#2 = 2}{
    \draw[draw=none, pattern=mydots, pattern color=purple] (0, -0.29*\s) rectangle (1*\s, -0.22*\s);
    }{}}

    {\ifthenelse{#3 = 2}{
    \draw[draw=none, pattern=mydots, pattern color=purple] (-0.29*\s, 0) rectangle (-0.21*\s, 1*\s);
    }{}}

    \ifthenelse{#2 = 1}{\draw[gray, thick] (0, -0.25*\s) -- (1*\s, -0.25*\s);}{\ifthenelse{#2 = 2}{\draw[blue, thick] (0, -0.29*\s) -- (1*\s, -0.29*\s);}{}}
    \draw[gray, thick] (0, -0.25*\s) -- (0, 0);
    \draw[gray, thick] (1*\s, -0.25*\s) -- (1*\s, 0);
    
    \ifthenelse{#3 = 1}{\draw[gray, thick] (-0.25*\s, 0) -- (-0.25*\s, 1*\s);}{\ifthenelse{#3 = 2}{\draw[red, thick] (-0.29*\s, 0) -- (-0.29*\s, 1*\s);}{}}
    \draw[gray, thick] (-0.25*\s, 0) -- (0, 0);
    \draw[gray, thick] (-0.25*\s, 1*\s) -- (0, 1*\s);

    \ifthenelse{#4 = 1}{\draw[gray, thick] (1.25*\s, 0) -- (1.25*\s, 1*\s);}{\ifthenelse{#4 = 2}{\draw[blue, thick] (1.29*\s, 0) -- (1.29*\s, 1*\s);}{}}
    \draw[gray, thick] (1.25*\s, 0) -- (1*\s, 0);
    \draw[gray, thick] (1.25*\s, 1*\s) -- (1*\s, 1*\s);

    \ifthenelse{#5 = 1}{\draw[gray, thick] (0, 1.25*\s) -- (1*\s, 1.25*\s);}{\ifthenelse{#5 = 2}{\draw[red, thick] (0, 1.29*\s) -- (1*\s, 1.29*\s);}{}}
    \draw[gray, thick] (0, 1*\s) -- (0, 1.25*\s);
    \draw[gray, thick] (1*\s, 1*\s) -- (1*\s, 1.25*\s);
    
    \draw[red!10, fill = red!10] (0, 0) rectangle (1*\s, 1*\s);
    
    \ifthenelse{#6 = 1}{\node at (0.5*\s, 0.5*\s) {$\Omega^c_{#1},\nu{=}\mu_{#1}$};}{\ifthenelse{#6 = 2}{\node at (0.5*\s, 0.5*\s) {$\nu{=}\hmu$};}}{\ifthenelse{#6 = 3}{\node at (0.5*\s, 0.5*\s) {$\hOmega^c$}; \node at (-0.1*\s, 0.5*\s) {$\hOmega^l$}; \node at (1.15*\s, 0.5*\s) {$\hOmega^r$}; \node at (0.5*\s, -0.12*\s) {$\hOmega^b$}; \node at (0.5*\s, 1.12*\s) {$\hOmega^t$};}{}}

}

\newcommand{\pateraWing}[4]{
    \ifthenelse{#1 = 1}{\draw[gray!5, fill = gray!5] (0, -0.25*\s) rectangle (1*\s, 0);}{\ifthenelse{#1 = 2}{\draw[blue!5, fill = blue!5] (0, -0.25*\s) rectangle (1*\s, 0);}}{\draw[gray!5, fill = gray!5] (0, -0.20*\s) rectangle (1*\s, 0);} 
    \ifthenelse{#2 = 1}{\draw[gray!5, fill = gray!5] (0, 1*\s) rectangle (1*\s, 1.25*\s);}{\ifthenelse{#2 = 2}{\draw[blue!5, fill = blue!5] (0, 1*\s) rectangle (1*\s, 1.25*\s);}}{\draw[gray!5, fill = gray!5] (0, 1*\s) rectangle (1*\s, 1.20*\s);} 
    \ifthenelse{#3 = 1}{\draw[gray!5, fill = gray!5] (-0.25*\s, 0) rectangle (0, 1*\s);}{\ifthenelse{#3 = 2}{\draw[yellow!5, fill = yellow!5] (-0.25*\s, 0) rectangle (0, 1*\s);}}{\draw[gray!5, fill = gray!5] (-0.20*\s, 0) rectangle (0, 1*\s);} 
    \ifthenelse{#4 = 1}{\draw[gray!5, fill = gray!5] (1*\s, 0) rectangle (1.25*\s, 1*\s);}{\ifthenelse{#4 = 2}{\draw[yellow!5, fill = yellow!5] (1*\s, 0) rectangle (1.25*\s, 1*\s);}}{\draw[gray!5, fill = gray!5] (1*\s, 0) rectangle (1.20*\s, 1*\s);} 

}

\newcommand{\all}{\text{\textit{:}}}

\graphicspath{{./figures/}}

\usepackage{a4wide}

\begin{document}

\begin{center}
\begin{Large}
\textbf{PGD-based local surrogate models via overlapping domain decomposition: a computational comparison}
\end{Large}

\bigskip

Marco Discacciati$^1$, Ben J. Evans$^1$, Matteo Giacomini$^{2,3}$

\medskip

\begin{small}
${}^1$ Department of Mathematical Sciences, Loughborough University, Epinal Way, LE11~3TU, Loughborough, United Kingdom. m.discacciati@lboro.ac.uk, b.j.evans@lboro.ac.uk

${}^2$ Laboratori de C\`alcul Numeric (LaC\`aN), E.T.S. de Ingenier\'ia de Caminos, Canales y Puertos, Universitat Polit\`ecnica de Catalunya - BarcelonaTech (UPC), Barcelona, Spain.
               
${}^3$ Centre Internacional de M\`etodes Num\`erics en Enginyeria (CIMNE), Barcelona, Spain. matteo.giacomini@upc.edu
\end{small}

\end{center}

\medskip

\begin{abstract}
An efficient strategy to construct physics-based local surrogate models for parametric linear elliptic problems is presented. The method relies on proper generalized decomposition (PGD) to reduce the dimensionality of the problem and on an overlapping domain decomposition (DD) strategy to decouple the spatial degrees of freedom.
In the offline phase, the local surrogate model is computed in a non-intrusive way, exploiting the linearity of the operator and imposing arbitrary Dirichlet conditions, independently at each node of the interface, by means of the traces of the finite element functions employed for the discretization inside the subdomain.
This leads to parametric subproblems with reduced dimensionality, significantly decreasing the complexity of the involved computations and achieving speed-ups up to 100 times with respect to a previously proposed DD-PGD algorithm that required clustering the interface nodes.
A fully algebraic alternating Schwarz method is then formulated to couple the subdomains in the online phase, leveraging the real-time (less than half a second) evaluation capabilities of the computed local surrogate models, that do not require the solution of any additional low-dimensional problems.
A computational comparison of different PGD-based local surrogate models is presented using a set of numerical benchmarks to showcase the superior performance of the proposed methodology, both in the offline and in the online phase.
\end{abstract}



\section{Introduction}
The construction of local reduced order models (ROMs) by means of domain decomposition (DD) approaches has gained increasing attention given the challenges of computing surrogate models of large-scale, possibly multi-physics, systems~\cite{Buhr:2021,Klawonn-HKLW-21}.
Both non-overlapping and overlapping DD methods have been studied in the literature, coupled with different ROM strategies, such as reduced basis (RB)~\cite{Maday:2002:JSC,Maday:2004:SISC,Iapichino:2012:CMAME,Patera-HKP-13,Iapichino:2016:CMA}, proper orthogonal decomposition (POD)~\cite{Baiges-BCI-13,Rozza-PNTBR-23,Iollo-IST-23}, and proper generalized decomposition (PGD)~\cite{Nazeer:2014:CM,Huerta:2017:IJNME,Discacciati:2024:CMAME}.

In order for these techniques to be suitable to treat realistic problems, three key aspects need to be fulfilled: (i) non-intrusiveness, (ii) physical interpretability, and (iii) efficiency.

Whilst projection-based ROMs encapsulate the physical information of the problem under analysis and can achieve significant dimensionality reduction and efficient performance~\cite{Buffoni:2009:CF,Corigliano-CDM-13}, they rely on the solution of a low-dimensional problem during the online phase to evaluate the surrogate model for a new set of parameters.
Such an intrusive implementation represents a potential bottleneck when commercial or industrial solvers are employed. To circumvent this issue, data-driven solutions constructing local functional approximations (e.g., with radial basis or Gaussian process regression) starting from the POD basis have been proposed in~\cite{Pain-XFHNP-19,Pain-XHFMHBANP-19}.
Nonetheless, data-driven surrogate models are known to suffer from limited physical interpretability, hence the need to introduce suitably-defined loss functions to include information on the underlying physics during training~\cite{Li-LTWL-19,Karniadakis-KZK-21,NissenMeyer-MMN-23}.
The main drawback of these approaches is represented by the large amount of data required for training and, consequently, by the high computational cost of constructing the surrogate models.

PGD~\cite{Chinesta-AMCK-06,Chinesta:2014} provides a competitive alternative to devise non-intrusive, physics-based surrogate models allowing efficient evaluations in real time.
This is achieved by means of an \textit{a priori} ROM framework in which the surrogate model is constructed during the offline phase as a rank-one approximation, the \textit{PGD expansion}, without the need to previously compute any snapshot of the problem.
Interested readers are referred to~\cite{Giacomini-GBSH-21} for a detailed comparison of \textit{a priori} and \textit{a posteriori} PGD strategies.
In particular, non-intrusive implementations have been successfully devised for PGD-ROMs using SAMCEF~\cite{Ladeveze-CNLB-16}, Abaqus~\cite{Zou-ZCDA-18}, OpenFOAM~\cite{Tsiolakis-TGSOH-20,Tsiolakis-TGSOH-22}, VPS/Pamcrash~\cite{Rocas-RGZLD-21}, and MSC-Nastran~\cite{Cavaliere-CZSLD-21,Cavaliere-CZSLD-22a}.
Moreover, the construction of the PGD expansion is driven by the minimization of a residual functional that measures the discrepancy between the rank-one approximation and the high-dimensional solution of the parametric partial differential equation (PDE) under analysis, thus naturally fulfilling the underlying physics.
Finally, the online phase only requires interpolation in the parametric space to retrieve the evaluation of the surrogate model for a new set of parameters, achieving efficient real-time performance by avoiding the solution of any additional problem.

The present work builds upon the approach coupling PGD-based surrogate models and the overlapping Schwarz method proposed in~\cite{Discacciati:2024:CMAME}.
We devise a novel PGD-ROM methodology significantly improving the performance of the original DD-PGD algorithm, and we perform a detailed computational comparison of the resulting local surrogate models.
Exploiting the linearity of the parametric PDE selected as model problem, the method defines local subproblems featuring arbitrary Dirichlet conditions at the interface.
The surrogate model with reduced dimensionality is thus obtained by constructing local ROMs with unitary boundary conditions at each node of the interface, using the traces of the finite element functions employed for the discretization within the subdomain.
This leads to a set of local subproblems featuring the same spatial and parametric dimensions as the original parametric PDE, which are solved using the non-intrusive, fully-algebraic framework of the Encapsulated PGD Algebraic Toolbox~\cite{Diez:2020:ACME}.
Hence, the proposed methodology outperforms existing PGD-based local ROMs: it allows to handle significantly smaller (and easier to solve) parametric local subproblems compared to the ones obtained in~\cite{Discacciati:2024:CMAME} by clustering the interface nodes; it circumvents the need for an expensive representation of the trace of the solution at the interface via auxiliary basis functions required by~\cite{Huerta:2017:IJNME}; it avoids the introduction of Lagrange multipliers (and their separated representations) used in~\cite{Nazeer:2014:CM} to couple the local surrogate models in the overlapping region.
Finally, in the online phase, the parametric linear system arising from imposing the equality of the traces of the local PGD solutions at the interfaces is solved by means of a matrix-free Krylov method, such as GMRES~\cite{Saad:1986:SISSC}. This is performed without the solution of any additional problem, thus allowing for the coupling to be executed in real time.

The rest of this article is structured as follows.
Section~\ref{sc:Problem} introduces the two-domain formulation of the parametric PDE used as model problem in this work, discusses the definition of the parametric trace of the solution at the interface, and recalls the algebraic form of the overlapping Schwarz algorithm.
In Sect.~\ref{sc:offlineOnline}, the offline phase of the construction of the local surrogate model with reduced dimensionality is presented. Moreover, the online phase featuring the coupling algorithm based on the surrogate model is described.
A critical discussion on the novel local surrogate model with reduced dimensionality and its comparison with the original DD-PGD strategy from~\cite{Discacciati:2024:CMAME} is presented in Sect.~\ref{sc:comparison}.
Section~\ref{sc:numericalResults} presents a computational comparison of the two approaches via a set of benchmark tests to assess the superior performance of the PGD-based local surrogate model with reduced dimensionality both in the offline and in the online phase.
Finally, Sect.~\ref{sc:Conclusion} summarizes the results of this work.

\section{Problem setting and overlapping domain decomposition}
\label{sc:Problem}

In this section, we introduce the model problem employed for the computational study of the different strategies to construct local surrogate models using PGD.
For the sake of clarity, we consider a parametric Poisson equation with Dirichlet boundary conditions, and we decompose the computational domain into two overlapping subdomains only. Extensions to the case of more general elliptic problems with different boundary conditions and to multiple subdomains can be achieved as explained in~\cite{Discacciati:2024:CMAME,Discacciati-DEG-25a}. Moreover, numerical results in Sect.~\ref{sc:numericalResults} will address different problems beyond the Poisson equation, showcasing the generality of the discussed methodologies.

\subsection{Two-domain formulation of the parametric Poisson equation}
\label{sc:twoDomain}
Let $\Omega \subset \mathbb{R}^d$ ($d = 1, 2, 3$) be an open bounded domain with Lipschitz boundary $\partial\Omega$, and let $\bmu =(\mu_1,\ldots,\mu_{N_p}) \in \mathcal{P}$ be a tuple of $N_p \in \mathbb{N}$ problem parameters with $\mathcal{P} = \mathcal{I}_1 \times \dots \times \mathcal{I}_{N_p} \subset \mathbb{R}^{N_p}$ and $\mathcal{I}_k$ compact ($k=1,\ldots,N_p$). Moreover, for all $\bmu \in \mathcal{P}$, let $\nu(\bmu) >0$ be a positive parametric diffusion coefficient, and $f(\bmu) \in L^2(\Omega)$ and $g(\bmu) \in H^{1/2}(\partial\Omega)$ be a parametric force and boundary datum, respectively. Then, consider the well-posed linear elliptic parametric boundary value problem: for all $\bmu \in \mathcal{P}$, find $u(\bmu) \in H^1(\Omega)$ such that 
\begin{subequations}\label{eq:globalProblem}
	\begin{eqnarray}
		-\nabla\cdot(\nu(\bmu)\nabla u(\bmu)) = f(\bmu) && \text{ in } \Omega, \label{eq:globalProblemA}\\
		u(\bmu) = g(\bmu) && \text{ on } \partial\Omega.
	\end{eqnarray}
\end{subequations}

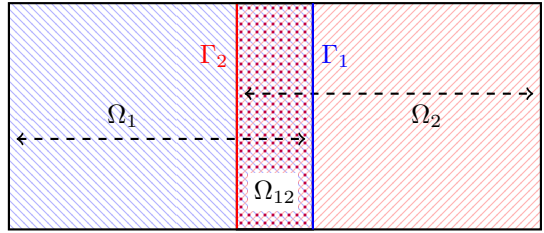
\begin{figure}[!htb]
	\centering
	\begin{tikzpicture}
        \draw[pattern=north west lines, pattern color=blue!30] (0, 0) rectangle (4, 3);
        \draw[pattern=north east lines, pattern color=red!30] (3, 0) rectangle (7, 3);        
		\draw[pattern=dots, pattern color=purple] (3, 0) rectangle (4, 3);
        \draw[fill=white, draw=none] (3.15, 0.25) rectangle (3.8, 0.75);
        \node at (3.5, 0.5) {$\Omega_{12}$};
				
        \node[blue] at (4.3, 2.3) {$\Gamma_1$};
		
        \node[red] at (2.7, 2.3) {$\Gamma_2$};
		
		\draw[black, thick] (3, 0) rectangle (7, 3);
        \draw[red, thick] (3, 0) -- (3, 3);
		\node[black] at (5.5, 1.5) {$\Omega_2$};
		\draw[<->, thick, black,dashed]  (3.1, 1.8) -- (6.9, 1.8);
		
		\draw[black, thick] (0, 0) rectangle (4, 3);
        \draw[blue, thick] (4, 0) -- (4, 3);
		\node[black] at (1.5, 1.5) {$\Omega_1$};
		\draw[<->, thick, black,dashed]  (0.1, 1.2) -- (3.9, 1.2);
		
	\end{tikzpicture}
\caption{Partition of the domain $\Omega$ into two overlapping subdomains $\Omega_1$ (light blue) and $\Omega_2$ (light red), with overlap $\Omega_{12}$ (dotted purple) and interfaces $\Gamma_1$ (blue) and $\Gamma_2$ (red).}\label{fig:partition}
\end{figure}

The domain $\Omega$ is split into two overlapping subdomains $\Omega_i \subset \Omega$ ($i = 1, 2$) such that $\Omega = \Omega_1 \cup \Omega_2$ and $\Omega_{12} = \Omega_1 \cap \Omega_2 \neq \emptyset$, and, for $i=1,2$, let $\Gamma_i = \partial\Omega_i \setminus \partial\Omega$ be the interfaces that are assumed not to intersect, that is, $dist(\Gamma_1,\Gamma_2)>0$ (see Fig.~\ref{fig:partition}). 

Considering this decomposition, problem~\eqref{eq:globalProblem} can be equivalently reformulated in the two-domain form:
for all $\bmu \in \mathcal{P}$, find $u_i(\bmu) \in H^1(\Omega_i)$ ($i = 1, 2$) such that
\begin{subequations}\label{eq:twoDomain}
	\begin{eqnarray}
		-\nabla\cdot(\nu(\bmu)\nabla u_i(\bmu)) = f_i(\bmu) && \text{ in } \Omega_i,\\
		u_i(\bmu) = g_i(\bmu) && \text{ on } \partial\Omega_i \cap \partial\Omega,\\
		u_1(\bmu) = u_2(\bmu) && \text{ on } \Gamma_1 \cup \Gamma_2, \label{eq:twoDomain_coupling}
	\end{eqnarray}
\end{subequations}
where $f_i(\bmu)$ and $g_i(\bmu)$ denote the restriction of $f(\bmu)$ to $\Omega_i$ and of $g(\bmu)$ to $\partial\Omega_i \cap \partial\Omega$, while~\eqref{eq:twoDomain_coupling} ensures the continuity of the local solutions $u_1(\bmu)$ and $u_2(\bmu)$ across the interfaces $\Gamma_1$ and $\Gamma_2$. 

\begin{rem}
In non-overlapping DD, both continuity of the solution and of the fluxes needs to be enforced at the interface for the two-domain formulation to be equivalent to the original problem~\cite{Quarteroni:1999,Toselli:2005}. This is not the case for overlapping DD and a proof of the equivalence between~\eqref{eq:globalProblem} and~\eqref{eq:twoDomain} can be found in~\cite{Discacciati:2013:SICON}.
This choice is particularly appealing from a computational viewpoint because it requires end-users to only have access to a numerical solver allowing to impose arbitrary Dirichlet boundary conditions on the newly introduced interfaces.
\end{rem}

We can then split the local problems~\eqref{eq:twoDomain} into two contributions.
\begin{enumerate}
\item
Local parametric problem depending on assigned data: find $u_i^f(\bmu) \in H^1(\Omega_i)$ such that
\begin{subequations}\label{eq:localProblemData}
	\begin{eqnarray}
		-\nabla\cdot(\nu(\bmu)\nabla u_i^f(\bmu)) = f_i(\bmu) && \text{ in } \Omega_i,\\
		u_i^f(\bmu) = g_i(\bmu) && \text{ on } \partial\Omega_i \cap \partial\Omega,\\
		u_i^f(\bmu) = \widetilde{g}_i(\bmu) && \text{ on } \Gamma_i,
	\end{eqnarray}
\end{subequations}
where $\widetilde{g}_i(\bmu) \in H^{1/2}(\Gamma_i)$ is a suitably defined trace function such that no discontinuity is introduced in the boundary condition on $\partial\Omega_i$ at the intersection between the interface and the external boundary.
\item
Local parametric problem depending on auxiliary interface data: find $u_i^\lambda(\bmu) \in H^1(\Omega_i)$ such that
\begin{subequations}\label{eq:localProblemLambda}
	\begin{eqnarray}
		-\nabla\cdot(\nu(\bmu)\nabla u_i^\lambda(\bmu)) = 0 && \text{ in } \Omega_i,\\
		u_i^\lambda(\bmu) = 0 && \text{ on } \partial\Omega_i \cap \partial\Omega,\\
		u_i^\lambda(\bmu) = \lambda_i(\bmu) && \text{ on } \Gamma_i, \label{eq:localProblemLambda3}
	\end{eqnarray}
\end{subequations}
where $\lambda_i \in H^{1/2}_{00}(\Gamma_i)$ must be chosen in such a way that
\begin{equation}\label{eq:funcSplitting}
u_i(\bmu) = u_i^f(\bmu) + u_i^\lambda(\bmu) \quad \text{for all } \bmu \in \mathcal{P}, \quad i=1,2.
\end{equation}
\end{enumerate}

\subsection{Definition of the parametric traces and finite element discretization}
\label{sc:Traces}
Problems~\eqref{eq:localProblemData} and~\eqref{eq:localProblemLambda} are now discretized by means of a user-selected solver. For the purpose of this work, a continuous Galerkin finite element method (FEM) is employed. For $i=1,2$, we consider two regular computational finite element meshes $\mathcal{T}_i$ in $\Omega_i$. For simplicity of presentation, we assume that the two meshes coincide in $\Omega_{12}$ and are conforming with the interfaces $\Gamma_1$ and $\Gamma_2$. 

To introduce the discrete space for the finite element discretization, we consider the space
\begin{equation*}
X_i^r = \{ v \in C^0(\overline{\Omega}_i) \, : \, v_{\vert K} \in \mathbb{P}_r \text{ or } v_{\vert K} \in \mathbb{Q}_r \text{ for all } K \in \mathcal{T}_i \} ,
\end{equation*}
where $\mathbb{P}_r$ denotes the space of polynomial functions of degree $\leq r$ for simplices, whereas $\mathbb{Q}_r$ is the tensor product polynomial space on quadrilateral/hexahedral elements, featuring polynomials of degree up to $r$ in each direction. Hence, the finite element space is defined as
\begin{equation*}
V_i^h = H^1_0(\Omega_i) \cap X_i^r \,.
\end{equation*}

Let $\{\varphi_i^j(\bx)\}_{j=1,\ldots,N_{\Omega_i}}$ be the finite element basis functions of $X_i^r$, and let $\eta_i^j (\bx) = \varphi_i^j(\bx)_{\vert\Gamma_i}$ be the restrictions of the finite element basis functions on the interface $\Gamma_i$. Note that the set of non-null trace functions
\begin{equation}\label{eq:femBasisTrace}
    \{ \eta_i^j(\bx)\}_{j=1,\ldots,N_{\Gamma_i}} ,
\end{equation}
with $N_{\Gamma_i}$ being the number of degrees of freedom (DOFs) on $\Gamma_i$, forms a partition of unity on $\Gamma_i$ and a basis for the discrete space of traces
\begin{equation*}
Y_i^h = \{ \lambda \in C^0(\overline{\Gamma}_i) \, : \, \lambda = 0 \text{ on } \overline{\Gamma}_i \cap \overline{\partial \Omega} \text{ and } \exists v \in X_i^r \text{ s.t. } v = \lambda \text{ on } \Gamma_i \}.
\end{equation*}
At the discrete level, we can define $g_{\Omega_i}^h (\bmu) \in H^1(\Omega_i) \cap X_i^r$ to be a continuous extension of the finite element interpolant of the boundary data in~\eqref{eq:localProblemData}. Similarly, for $i=1,2$, we can approximate the trace function at the interface in \eqref{eq:localProblemLambda3} as
\begin{equation}\label{eq:lambdaRepresentation}
\lambda_i (\bmu) \approx \lambda_i^h (\bmu) = \lambda_i^h (\bx;\bmu) = \sum_{j=1}^{N_{\Gamma_i}} \Lambda_i^j (\bmu) \, \eta_i^j (\bx) \; \text{ on }\Gamma_i,
\end{equation}
whereas its discrete extension within the subdomain $\Omega_i$ can be constructed as
\begin{equation}\label{eq:lambdaExtension}
\lambda_{\Omega_i}^h (\bmu) = \lambda_{\Omega_i}^h (\bx;\bmu) = \sum_{j=1}^{N_{\Gamma_i}} \Lambda_i^j (\bmu) \, \varphi_i^j (\bx) \; \text{ in } \Omega_i ,
\end{equation}
where $\{ \Lambda_i^j (\bmu) \}_{j=1, \ldots,N_{\Gamma_i}}$ denote the nodal values of $\lambda_i(\bmu)$ at the $N_{\Gamma_i}$ DOFs on $\Gamma_i$. Note that, by construction, $\lambda_{\Omega_i}^h (\bmu) = \lambda_i^h (\bmu)$ on $\Gamma_i$.

Hence, we can express the Galerkin approximation of $u_i^f (\bmu)$ and $u_i^\lambda(\bmu)$ in~\eqref{eq:funcSplitting} as the superposition of the extension of the Dirichlet data $g_{\Omega_i}^h (\bmu)$ (respectively, $\lambda_{\Omega_i}^h (\bmu)$) and the new finite element function $u_{i,h}^{0,f}(\bmu)$ (respectively, $u_{i,h}^{0,\lambda}(\bmu)$), belonging to the space $V_i^h$, namely,
\begin{equation}\label{eq:galerkinApproximations}
u_{i,h}^f (\bmu) = u_{i,h}^{0,f} (\bmu) + g_{\Omega_i}^h (\bmu)
\quad \text{ and } \quad
u_{i,h}^\lambda (\bmu) = u_{i,h}^{0,\lambda} (\bmu) + \lambda_{\Omega_i}^h (\bmu)\,.
\end{equation}
It follows that the Galerkin approximation of $u_i(\bmu)$ in~\eqref{eq:twoDomain} becomes
\begin{equation}\label{eq:uLambdaData}
u_{i,h}(\bmu) = u_{i,h}^0 (\bmu) + \lambda_{\Omega_i}^h (\bmu) + g_{\Omega_i}^h (\bmu) ,
\end{equation}
where
\begin{equation}\label{eq:solSuperposition}
u_{i,h}^0(\bmu) = u_{i,h}^{0,f}(\bmu) + u_{i,h}^{0,\lambda}(\bmu).
\end{equation}

Exploiting~\eqref{eq:uLambdaData},  the coupling conditions~\eqref{eq:twoDomain_coupling} can be reformulated at the discrete level as
\begin{subequations}\label{eq:discreteInterfaceConditions}
\begin{eqnarray}
\lambda_1^h(\bmu) - \lambda_{\Omega_2}^h (\bmu)_{\vert\Gamma_1} - u_{2,h}^0 (\bmu)_{\vert\Gamma_1} = (-g_{\Omega_1}^h(\bmu) + g_{\Omega_2}^h(\bmu))_{\vert\Gamma_1} && \text{ on } \Gamma_1, \label{eq:discreteInterfaceConditions1} \\
\lambda_2^h(\bmu) - \lambda_{\Omega_1}^h (\bmu)_{\vert\Gamma_2} - u_{1,h}^0 (\bmu)_{\vert\Gamma_2} = (\phantom{-}g_{\Omega_1}^h (\bmu) - g_{\Omega_2}^h(\bmu))_{\vert\Gamma_2} && \text{ on } \Gamma_2.
\end{eqnarray}
\end{subequations}
Therefore, the Galerkin finite element approximation of the two-domain formulation~\eqref{eq:twoDomain} becomes: for all $\bmu \in \mathcal{P}$, for $i=1,2$, find $u_{i,h}^0 (\bmu) \in V_i^h$ and $\lambda_i^h (\bmu) \in Y_i^h$ such that the local problems
\begin{equation}\label{eq:galerkinTwoDomain}
\mathcal{A}_i (u_{i,h}^0(\bmu) + \lambda_{\Omega_i}^h (\bmu), v_{i,h} ; \bmu) = \mathcal{F}_i(v_{i,h}; \bmu) - \mathcal{A}_i (g_{\Omega_i}^h (\bmu) , v_{i,h} ; \bmu) ,
\end{equation}
and the interface conditions~\eqref{eq:discreteInterfaceConditions} are satisfied for all $v_{i,h} \in V_i^h$, where, for all $u,v \in H^1(\Omega_i)$, the following variational forms are introduced
\begin{equation}\label{eq:defForms}
\mathcal{A}_i(u,v;\bmu) = \int_{\Omega_i} \nu(\bmu) \nabla u \cdot \nabla v \, d\bx
\quad \text{ and } \quad
\mathcal{F}_i(v; \bmu) = \int_{\Omega_i} f_i(\bmu) v \, d\bx .
\end{equation}

\subsection{Algebraic form of the alternating Schwarz method}
\label{sc:AlgebraicForm}
Starting from the discrete problem~\eqref{eq:galerkinTwoDomain} with interface conditions~\eqref{eq:discreteInterfaceConditions}, this section formulates an algebraic substructuring version of the classical alternating Schwarz method, following ideas from~\cite{Smith:1996}. Note that this algorithm is the basis for the online coupling phase of the local surrogate models presented in the following sections.

Let $\mat{A}_{\all\all}^i$ be the stiffness finite element matrix associated with the bilinear form $\mathcal{A}_i(\cdot,\cdot;\bmu)$ in $\Omega_i$. We replace the sub-index $\all$ by $I$ (respectively, $\Gamma_i$) to denote that only the rows/columns associated with the DOFs internal to $\Omega_i$ excluding $\Gamma_i$ (respectively, the DOFs on $\Gamma_i$) are retained, with the index for rows preceding the one for columns. Also, let $\mat{I}_{\Gamma_i\Gamma_i}$ be the identity matrix associated with the DOFs on $\Gamma_i$, and $\mat{R}_{\Omega_i\to\Gamma_j}$ ($i,j=1,2$, $i \not= j$) be the restriction matrix that, given any vector of nodal values in $\Omega_i$, returns the vector of nodal values on the interface $\Gamma_j$ inside $\Omega_i$. 
Moreover, $\mat{E}_{\Gamma_i\to\Omega_i} \bLambda_{\Gamma_i} \in X_i^r$, with $\mat{E}_{\Gamma_i\to\Omega_i} \bLambda_{\Gamma_i} = \mathbf{0}$ on $\partial \Omega_i \cap \partial\Omega$, is a suitable algebraic extension operator of the interface nodal values $\bLambda_{\Gamma_i}$ that corresponds, e.g., to~\eqref{eq:lambdaExtension} at the algebraic level.
Finally, let $\mathbf{u}_I^i$ denote the vector of the nodal values of the finite element function $u_{i,h}^0 (\bmu)$.

Then, with additional self-explanatory notation, the linear system associated with the local subproblems~\eqref{eq:galerkinTwoDomain} for $i=1,2$ and the corresponding interface conditions~\eqref{eq:discreteInterfaceConditions} becomes
\begin{equation}\label{eq:twoDomainAlgebraic}
 \begin{pmatrix}
 \mat{A}^1_{II} & \mat{0} & \mat{A}^1_{I\Gamma_1} & \mat{0} \\[2pt]
 \mat{0} & \mat{A}^2_{II} & \mat{0} & \mat{A}^2_{I\Gamma_2} \\[2pt]
 \mat{0} & -\mat{R}_{\Omega_2 \to \Gamma_1} & \mat{I}_{\Gamma_1\Gamma_1} & -\mat{R}_{\Omega_2 \to \Gamma_1} \mat{E}_{\Gamma_2 \to \Omega_2} \\[2pt]
 -\mat{R}_{\Omega_1\to\Gamma_2} & \mat{0} & -\mat{R}_{\Omega_1 \to \Gamma_2} \mat{E}_{\Gamma_1 \to \Omega_1} & \mat{I}_{\Gamma_2\Gamma_2}
 \end{pmatrix}
 \begin{pmatrix}
     \mathbf{u}^1_I \\[2pt]
     \mathbf{u}^2_I \\[2pt]
     \boldsymbol{\Lambda}_{\Gamma_1} \\[2pt]
     \boldsymbol{\Lambda}_{\Gamma_2}
 \end{pmatrix}
 =
 \begin{pmatrix}
     \mathbf{f}^1_I \\[2pt]
     \mathbf{f}^2_I \\[2pt]
     \mathbf{g}_{\Gamma_1} \\[2pt]
     \mathbf{g}_{\Gamma_2}
 \end{pmatrix}.
\end{equation}

Following~\cite{Discacciati:2024:CMAME}, we can reduce~\eqref{eq:twoDomainAlgebraic} to the equivalent interface system
\begin{equation}\label{eq:twoDomainInterfaceAlgebraic}
 \Sigma
 \begin{pmatrix}
     \bLambda_{\Gamma_1} \\[2pt]
     \bLambda_{\Gamma_2}
 \end{pmatrix} \\
 =
 \begin{pmatrix}
     \mathbf{g}_{\Gamma_1} + \mat{R}_{\Omega_2 \to \Gamma_1} (\mat{A}^2_{II})^{-1} \mathbf{f}^2_I \\[2pt]
     \mathbf{g}_{\Gamma_2} + \mat{R}_{\Omega_1 \to \Gamma_2} (\mat{A}^1_{II})^{-1} \mathbf{f}^1_I
 \end{pmatrix},
\end{equation}
with
\begin{equation*}
\Sigma = 
 \begin{pmatrix}
 \mat{I}_{\Gamma_1\Gamma_1} & -\mat{R}_{\Omega_2 \to \Gamma_1} (\mat{E}_{\Gamma_2 \to \Omega_2} + (\mat{A}^2_{II})^{-1} (-\mat{A}^2_{I\Gamma_2}) ) \\
 -\mat{R}_{\Omega_1 \to \Gamma_2} (\mat{E}_{\Gamma_1 \to \Omega_1} + (\mat{A}^1_{II})^{-1} (-\mat{A}^1_{I\Gamma_1}) ) & \mat{I}_{\Gamma_2\Gamma_2}
 \end{pmatrix}.
\end{equation*}
System~\eqref{eq:twoDomainInterfaceAlgebraic} can be solved using a suitable matrix-free Krylov method (e.g., GMRES). 
Note that the expensive part of this iterative algorithm is the solution of the local problems in $\Omega_i$ ($i=1,2$), that is, the computation of
\begin{equation}\label{eq:subEqLoc}
\mathbf{u}^i_I = (\mat{A}^i_{II})^{-1} \mathbf{f}^i_I + (\mat{A}^i_{II})^{-1} (-\mat{A}^i_{I\Gamma_i}) \bLambda_{\Gamma_i} .
\end{equation}
Indeed, computing
\begin{subequations}\label{eq:localProblemsAlgebraic}
\begin{eqnarray}
&& (\mat{A}^i_{II})^{-1} \mathbf{f}^i_I, \label{eq:localProblemsAlgebraicData} \\
&& (\mat{A}^i_{II})^{-1} (-\mat{A}^i_{I\Gamma_i}) \bLambda_{\Gamma_i}\label{eq:localProblemsAlgebraicLambda}
\end{eqnarray}
\end{subequations}
can become especially demanding when a new set of parameters $\bmu$ must be considered because, in general, the matrices $\mat{A}_{II}^i$ and $\mat{A}_{I\Gamma_i}^i$ can depend on $\bmu$, so the entire procedure - including the matrix assembly - must be re-executed from scratch.

To reduce the computational cost of the coupling procedure~\eqref{eq:twoDomainInterfaceAlgebraic}, following~\cite{Discacciati:2024:CMAME}, we split the algorithm into an offline  and an online phase to efficiently handle the presence of the problem parameters $\bmu \in \mathcal{P}$. More precisely, in the offline phase, the local problems~\eqref{eq:localProblemsAlgebraic} are solved by PGD to construct physics-based local surrogate models. These will be used in the online phase to reformulate~\eqref{eq:twoDomainInterfaceAlgebraic} in such a way that, at each GMRES iteration, the solution of the local problems is replaced by the evaluation of the precomputed local surrogate models, thus reducing the overall computational cost.
The construction of the PGD-based local surrogate models is discussed in Sect.~\ref{sc:offlinePhase}, while the online phase is presented in Sect.~\ref{sc:onlinePhase}.

\section{Accelerating domain decomposition via surrogate models}
\label{sc:offlineOnline}

In this section, we present a surrogate-based strategy for real-time overlapping domain decomposition, accelerated by means of the precomputation of ROMs for the local subproblems.

\subsection{Offline phase: construction of PGD-based local surrogate models}
\label{sc:offlinePhase}
First, the surrogate models for the local subproblems~\eqref{eq:galerkinTwoDomain} are computed.
Recall that the parametric solution of the local problem can be decomposed as in \eqref{eq:solSuperposition}, with $u_{i,h}^{0,f}(\bmu)$ depending on the data of the problem and $u_{i,h}^{0,\lambda}(\bmu)$ on the trace of the solution at the interface.

\subsubsection{Separated form of the data problem}
\label{sc:DataSubPb}
The component $u_{i,h}^{0,f}(\bmu)$ of the solution~\eqref{eq:solSuperposition} is associated with the algebraic system~\eqref{eq:localProblemsAlgebraicData}.
This is obtained from the Galerkin discretization of the following parametric problem: for all $\bmu \in \mathcal{P}$, find $u_{i,h}^{0,f} (\bmu) \in V_i^h$ such that
\begin{equation}\label{eq:galerkinForce}
\mathcal{A}_i (u_{i,h}^{0,f}(\bmu), v_{i,h} ; \bmu ) = \mathcal{F}_i (v_{i,h}) - \mathcal{A}_i (g_{\Omega_i}^h (\bmu), v_{i,h} ; \bmu) \qquad \forall v_{i,h} \in V_i^h.
\end{equation}

Following the standard approach in PGD~\cite{Chinesta:2014}, we assume that all data are given in separated form as the sum of products of functions that depend either on the spatial coordinate $\bx$ or on the parameters $\bmu$. Therefore, we write
\begin{equation*}
\nu (\bmu) = \sum_{\ell = 1}^{N_\nu} \xi_\nu^\ell (\bmu) b_\nu^\ell (\bx), \quad
f_i (\bmu) = \sum_{\ell = 1}^{N_f}   \xi_{i,f}^\ell (\bmu) b_{i,f}^\ell (\bx), \quad
g_{\Omega_i}^h (\bmu) = \sum_{\ell=1}^{N_D}  \xi_{i,D}^\ell(\bmu) b_{i,D}^\ell(\bx) ,
\end{equation*}
with the parametric modes expressed as the product of scalar functions, each depending on one parameter $\mu_k, \ k=1,\ldots,N_p$, e.g.,
\begin{equation*}
\xi_\nu^\ell (\bmu) = \prod_{k=1}^{N_p} \xi_\nu^{\ell,k} (\mu_k) \, .
\end{equation*}
The contribution of the Dirichlet boundary condition $g_{\Omega_i}^h (\bmu)$ is handled by introducing \textit{ad-hoc}, sufficiently smooth modes as usually done in the PGD context~\cite{Chinesta:2014}. 

Moreover, the solution $u^{0,f}_{i,h}(\bmu)$ of~\eqref{eq:galerkinForce} is approximated in separated form as
\begin{equation}\label{eq:SepSolnU}
    u^{0,f}_{i,h}(\bmu) \approx u^{0,f}_{i,\pgd}(\bmu) = \displaystyle\sum_{m = 1}^{M_i^f} {U}_{i}^m(\bx)\,{\phi}_{i}^m(\bmu) \,,
\end{equation}
where ${U}_{i}^m(\bx)$ is the $m$th spatial mode that is discretized by the Galerkin finite element method, while ${\phi}_{i}^m(\bmu)$ is the corresponding parametric mode that is discretized by pointwise collocation.

Assuming an affine parameter dependence for the variational forms in~\eqref{eq:defForms}, for all $u,v \in H^1(\Omega_i)$ and $\bmu \in \mathcal{P}$, we define
\begin{equation}\label{eq:PGDoperator}
 \mathcal{A}_i^{\pgd}(u, v; \bmu) =
  \sum_{\ell=1}^{N_{\nu}} \xi_{\nu}^\ell(\bmu) \int_{\Omega_i} b_{\nu}^\ell(\bx)\, \nabla u \cdot \nabla v \, d\bx\,,
  \qquad
 \mathcal{F}_i^{\pgd} (v; \bmu) = \displaystyle \sum_{\ell=1}^{N_f} \xi_{i, f}^\ell(\bmu) \int_{\Omega_i} b_{i, f}^\ell(\bx) v \, d\bx  \,.
\end{equation}
Then, the PGD approximation ${u}^{0,f}_{i,\pgd}(\bmu)$ is computed by solving the parametric problem
\begin{equation}\label{eq:problemLocalPGDData}
 \mathcal{A}_i^{\pgd}({u}^{0,f}_{i,\pgd}(\bmu), v_{i,h}; \bmu) = \mathcal{F}_i^{\pgd}(v_{i,h};\bmu) - \mathcal{A}_i^{\pgd} \left(\sum_{\ell=1}^{N_D}  \xi_{i,D}^\ell(\bmu) b_{i,D}^\ell(\bx), v_{i,h} ; \bmu\right)\,,\end{equation}
for all $v_{i,h} \in V_i^h$ and for all $\bmu \in \mathcal{P}$, using the greedy strategy based on the alternating direction algorithm described in~\cite{Diez:2020:ACME}.

\subsubsection{Separated form of the interface problem with reduced dimensionality}
\label{sc:OfflineRedDim}
The approximation of the second component of~\eqref{eq:solSuperposition}, that is, $u_{i,h}^{0,\lambda}(\bmu)$, follows the same rationale presented in Sect.~\ref{sc:DataSubPb}.
First, we notice that system~\eqref{eq:localProblemsAlgebraicLambda} corresponds to the following Galerkin problem: for all $\bmu \in \mathcal{P}$, find $u_{i,h}^{0,\lambda}(\bmu) \in V_i^h$ such that
\begin{equation}\label{eq:galerkinLambda}
\mathcal{A}_i(u_{i,h}^{0,\lambda}(\bmu), v_{i,h} ; \bmu ) = - \mathcal{A}_i(\lambda_{\Omega_i}^h (\bmu), v_{i,h} ; \bmu) \qquad \forall v_{i,h} \in V_i^h,
\end{equation}
where $\lambda_{\Omega_i}^h(\bmu)$ is expressed in separated form as in~\eqref{eq:lambdaExtension}, with the coefficients $\{ \Lambda_i^j(\bmu) \}_{j=1,\ldots,N_{\Gamma_i}}$ being the elements of the vector $\bLambda_{\Gamma_i}$.

It is worth noticing that, in order to fully characterize the solution $u_{i,h}^{0,\lambda}(\bmu)$, the $N_{\Gamma_i}$ coefficients $\{ \Lambda_i^j(\bmu) \}_{j=1,\ldots,N_{\Gamma_i}}$ need to be determined.
Following classical ROM approaches, such coefficients can be considered as additional independent variables. Nonetheless, this quickly leads to extremely high-dimensional problems (of the order of several tens of parameters) which are computationally unaffordable, even for settings in two spatial dimensions.

In this work, we propose to significantly reduce the dimensionality of the local surrogate models by exploiting superposition and linearity of the parametric PDE under analysis.
More precisely, consider the basis functions \eqref{eq:femBasisTrace} on the interface $\Gamma_i$ and their associated finite element basis functions $\{ \varphi_i^j (\bx)\}_{j=1,\ldots,N_{\Gamma_i}}$ in $\Omega_i$.
Using the expansions~\eqref{eq:lambdaRepresentation} and~\eqref{eq:lambdaExtension}, and exploiting the linearity of the problem, we can define a set of $N_{\Gamma_i}$ new functions $u_{i,h}^j(\bmu) \in V_i^h$ such that, for all $\bmu \in \mathcal{P}$, it holds
\begin{equation}\label{eq:unitaryProblems}
\mathcal{A}_i(u_{i,h}^j(\bmu),v_{i,h};\bmu) = - \mathcal{A}_i(\varphi_i^j (\bx), v_{i,h} ; \bmu) \qquad \forall v_{i,h} \in V_i^h,
\end{equation}
and 
\begin{equation}\label{eq:ulambdaExpansion}
u_{i,h}^{0,\lambda} (\bmu) = \sum_{j=1}^{N_{\Gamma_i}} \Lambda_i^j (\bmu) \, u_{i,h}^j (\bmu) .
\end{equation}

Following the same rationale as in~\eqref{eq:SepSolnU}, for $j=1,\ldots,N_{\Gamma_i}$, we introduce the separated representation of the solution $u_{i,h}^j(\bmu)$ defined as
\begin{equation}\label{eq:SepSolnUlambda}
    u_{i,h}^j(\bmu) \approx u_{i,\pgd}^j(\bmu) = \displaystyle\sum_{m = 1}^{M_i^j} {U}_{i,j}^m(\bx)\,{\phi}_{i,j}^m(\bmu) \,,
\end{equation}
where ${U}_{i,j}^m(\bx)$ and ${\phi}_{i, j}^m(\bmu)$ denote the $m$th spatial and parametric modes in the PGD expansion.
The local surrogate model $u^{j}_{i,\pgd}(\bmu)$ is then computed by solving the local parametric problem approximating~\eqref{eq:unitaryProblems}, namely,
\begin{equation}\label{eq:problemLocalPGDUnity}
 \mathcal{A}_i^{\pgd}(u^{j}_{i,\pgd}(\bmu), v_{i,h}; \bmu) = - \mathcal{A}_i^{\pgd} (\varphi_i^j(\bx), v_{i,h} ; \bmu)\,, \qquad \forall v_{i,h} \in V_i^h, \; \bmu \in \mathcal{P} ,
\end{equation}
with the separated PGD operators introduced in~\eqref{eq:PGDoperator}.

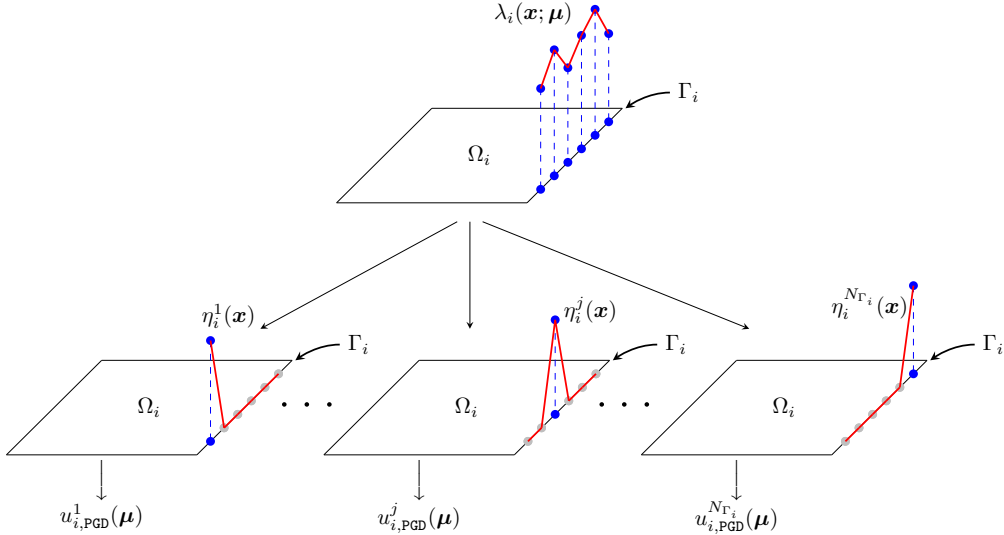
\begin{figure}[!ht]
    \centering
    \resizebox{0.9\textwidth}{!}{
    \begin{tikzpicture}
			\newcommand{\intNodes}[3]{
				\draw[dashed, #3] (\longSide + #1 * \delta / 7, \ang * #1 * \delta / 7 + #2) -- (\longSide + #1 * \delta / 7, #1 * \ang * \delta / 7 );
				\fill[#3] (\longSide + #1 * \delta / 7, \ang * #1 * \delta / 7 + #2) circle (0.07 cm);
				\fill[#3] (\longSide + #1 * \delta / 7, #1 * \ang * \delta / 7 ) circle (0.07 cm);
			}
			
			\newcommand{\intFun}[3]{
				\draw[thick, red] (\longSide + #1 * \delta / 7, #1 * \ang * \delta / 7 + #2 ) -- (\longSide + #1 * \delta / 7 + \delta / 7, #1 * \ang * \delta / 7 +  \ang * \delta / 7 + #3);
			}
			
			\def\ang{1}
			\def\delta{1.5}
			\def\longSide{3}
			
			\begin{scope}[shift={(0.2, 0)}]
				\draw (0, 0) -- (\delta, \delta*\ang);
				\draw (0, 0) -- (\longSide, 0);
				\draw(\delta, \delta*\ang) -- (\delta + \longSide, \delta*\ang);
				\draw(\longSide, 0) -- (\delta + \longSide, \delta*\ang);
				
				\node at (\longSide / 2 + \delta / 2, \delta * \ang / 2) {$\Omega_i$};
				
				\draw[-stealth, thick] (\longSide + \delta + 0.75, \delta * \ang + 0.25) node[anchor=west] {$\Gamma_i$} parabola (\longSide + \delta + 0.05, \delta * \ang);
				
				\intNodes{1}{1.6}{blue};
				\intNodes{2}{2}{blue};
				\intNodes{3}{1.5}{blue};	
				\intNodes{4}{1.8}{blue};
				\intNodes{5}{2}{blue};
				\intNodes{6}{1.4}{blue};
				
				\intFun{1}{1.6}{2};
				\intFun{2}{2}{1.5};
				\intFun{3}{1.5}{1.8};
				\intFun{4}{1.8}{2};
				\intFun{5}{2}{1.4};
			\end{scope}
			
			\begin{scope}[shift={(-5, -4)}]
				\draw (0, 0) -- (\delta, \delta*\ang);
				\draw (0, 0) -- (\longSide, 0);
				\draw(\delta, \delta*\ang) -- (\delta + \longSide, \delta*\ang);
				\draw(\longSide, 0) -- (\delta + \longSide, \delta*\ang);
				
				\node at (\longSide / 2 + \delta / 2, \delta * \ang / 2) {$\Omega_i$};
				
				\draw[-stealth, thick] (\longSide + \delta + 0.75, \delta * \ang + 0.25) node[anchor=west] {$\Gamma_i$} parabola (\longSide + \delta + 0.05, \delta * \ang);
				
				\intNodes{1}{1.6}{blue};
				\intNodes{2}{0}{lightgray};
				\intNodes{3}{0}{lightgray};
				\intNodes{4}{0}{lightgray};
				\intNodes{5}{0}{lightgray};
				\intNodes{6}{0}{lightgray};

				\intFun{1}{1.6}{0};
				\intFun{2}{0}{0};
				\intFun{3}{0}{0};
				\intFun{4}{0}{0};
				\intFun{5}{0}{0};
			\end{scope}
		
			\begin{scope}[shift={(0, -4)}]
				\draw (0, 0) -- (\delta, \delta*\ang);
				\draw (0, 0) -- (\longSide, 0);
				\draw(\delta, \delta*\ang) -- (\delta + \longSide, \delta*\ang);
				\draw(\longSide, 0) -- (\delta + \longSide, \delta*\ang);
				
				\node at (\longSide / 2 + \delta / 2, \delta * \ang / 2) {$\Omega_i$};
				
				\draw[-stealth, thick] (\longSide + \delta + 0.75, \delta * \ang + 0.25) node[anchor=west] {$\Gamma_i$} parabola (\longSide + \delta + 0.05, \delta * \ang);
				
				\intNodes{1}{0}{lightgray};
				\intNodes{2}{0}{lightgray};
				
				\intNodes{3}{1.5}{blue};
				
				\intNodes{4}{0}{lightgray};
				\intNodes{5}{0}{lightgray};
				\intNodes{6}{0}{lightgray};

				\intFun{1}{0}{0};
				\intFun{2}{0}{1.5};
				\intFun{3}{1.5}{0};
				\intFun{4}{0}{0};
				\intFun{5}{0}{0};
			\end{scope}
		
			\begin{scope}[shift={(5, -4)}]
				\draw (0, 0) -- (\delta, \delta*\ang);
				\draw (0, 0) -- (\longSide, 0);
				\draw(\delta, \delta*\ang) -- (\delta + \longSide, \delta*\ang);
				\draw(\longSide, 0) -- (\delta + \longSide, \delta*\ang);
				
				\node at (\longSide / 2 + \delta / 2, \delta * \ang / 2) {$\Omega_i$};
				
				\draw[-stealth, thick] (\longSide + \delta + 0.75, \delta * \ang + 0.25) node[anchor=west] {$\Gamma_i$} parabola (\longSide + \delta + 0.05, \delta * \ang);
				
				\intNodes{1}{0}{lightgray};
				\intNodes{2}{0}{lightgray};	
				\intNodes{3}{0}{lightgray};
				\intNodes{4}{0}{lightgray};
				\intNodes{5}{0}{lightgray};
				\intNodes{6}{1.4}{blue};				
				
				\intFun{1}{0}{0};
				\intFun{2}{0}{0};
				\intFun{3}{0}{0};
				\intFun{4}{0}{0};
				\intFun{5}{0}{1.4};
			\end{scope}
			
			\draw[- stealth] (2.1, -0.3) -- (-1 , -2);
			\draw[- stealth] (2.3, -0.3) -- (2.3 , -2);
			\draw[- stealth] (2.5, -0.3) -- (6.7 , -2);
			
			\node at (-0.2, -3.2) {\Huge \ldots};
			\node at (4.8, -3.2) {\Huge \ldots};
			
			\node at (-3.5, -4.4) {$\Big\downarrow$};
			\node at ( 1.5, -4.4) {$\Big\downarrow$};
			\node at ( 6.5, -4.4) {$\Big\downarrow$};
			
			\node at (-3.5, -5.0) {${u}_{i,\pgd}^1(\bmu)$};
			\node at ( 1.5, -5.0) {${u}_{i,\pgd}^j(\bmu)$};
			\node at ( 6.5, -5.0) {${u}_{i,\pgd}^{N_{\Gamma_i}}(\bmu)$};

                \node at ( 3.3,  3.0) {$\lambda_i(\bx;\bmu)$};
                \node at (-1.5, -1.8) {$\eta_i^1(\bx)$};
                \node at ( 4.2, -1.7) {$\eta_i^j(\bx)$};
                \node at ( 8.6, -1.6) {$\eta_i^{N_{\Gamma_i}}(\bx)$};
			
		\end{tikzpicture}
  }
    \caption{Partition of the interface nodes as a collection of single independent interface parameters.}
    \label{fig:lambdaSplitting}
\end{figure}

Recalling that $\eta_i^j (\bx) = \varphi_i^j(\bx)_{\vert\Gamma_i}$, Fig.~\ref{fig:lambdaSplitting} displays a schematic representation of the decomposition of the interface datum as a superposition of local spatial modes and the corresponding PGD expansions $\{ u_{i,\pgd}^j (\bmu) \}_{j=1,\ldots,N_{\Gamma_i}}$.
Note that the resulting problems~\eqref{eq:problemLocalPGDUnity} can be solved independently from one another and do not depend on the $N_{\Gamma_i}$ coefficients $\{ \Lambda_i^j(\bmu) \}_{j=1,\ldots,N_{\Gamma_i}}$.
This leads to local subproblems of the same spatial and parametric dimensions as the original parametric PDE~\eqref{eq:globalProblem}, circumventing the previously-mentioned issue of increasing dimensionality.

\begin{rem}
For $i=1,2$, if $\bLambda_{\Gamma_i} = (\Lambda_i^1 (\bmu), \ldots, \Lambda_i^{N_{\Gamma_i}} (\bmu) )^T$ are the solutions of system~\eqref{eq:twoDomainInterfaceAlgebraic}, it holds that $u_{i,h}^0 (\bmu) = u_{i,h}^{0,f} (\bmu) + u_{i,h}^{0,\lambda}(\bmu)$.    
\end{rem}

The function ${u}^{0,f}_{i,\pgd}(\bmu)$ obtained from~\eqref{eq:problemLocalPGDData} and the pairs $\{(u^{j}_{i,\pgd}(\bmu),\Lambda_i^j(\bmu))\}_{j = 1, \dots, N_{\Gamma_i}}$ consisting of the PGD expansion computed by solving~\eqref{eq:problemLocalPGDUnity} and the set of coefficients in the vector $\bLambda_{\Gamma_i}$ form the PGD-based local surrogate model for problem~\eqref{eq:galerkinTwoDomain}.
This provides a computationally inexpensive representation of the solution of the parametric problem~\eqref{eq:globalProblem} in subdomain $\Omega_i$.
Note that whilst ${u}^{0,f}_{i,\pgd}(\bmu)$ and $u^{j}_{i,\pgd}(\bmu)$ are fully characterized by the variational problems mentioned above, the coefficients $\{ \Lambda_i^j(\bmu) \}_{j=1,\ldots,N_{\Gamma_i}}$ are yet to be determined.
The corresponding procedure is described in Sect.~\ref{sc:onlinePhase}.

\subsection{Online phase: surrogate-based coupling procedure}
\label{sc:onlinePhase}
To reduce the computational cost of the coupling procedure~\eqref{eq:twoDomainInterfaceAlgebraic}, we utilize the PGD local surrogate models computed in the offline phase. 

In particular, recall that system~\eqref{eq:localProblemsAlgebraicData} corresponds to problem~\eqref{eq:galerkinForce} whose solution is approximated by the local surrogate model $u_{i,\pgd}^{0,f}(\bmu)$ computed in~\eqref{eq:problemLocalPGDData}.
We can thus replace the solution of~\eqref{eq:localProblemsAlgebraicData} by the evaluation of the surrogate model $u_{i,\pgd}^{0,f}(\bmu)$ for a suitable value, say, $\bar{\bmu} \in \mathcal{P}$, of the parameters, i.e., $u_{i,\pgd}^{0,f}(\bar{\bmu})$.

Similarly, we can exploit the correspondence between~\eqref{eq:localProblemsAlgebraicLambda},~\eqref{eq:galerkinLambda} and the set of problems~\eqref{eq:problemLocalPGDUnity}. 
More precisely, using~\eqref{eq:ulambdaExpansion} and~\eqref{eq:SepSolnUlambda} we can replace, for any vector $\bLambda_{\Gamma_i}$, the solution of~\eqref{eq:localProblemsAlgebraicLambda} by the local PGD operator
\begin{equation}\label{eq:poissonPGD}
\mathcal{P}_i^\pgd: \bLambda_{\Gamma_i} \to \sum_{j=1}^{N_{\Gamma_i}} \Lambda_i^j(\bmu) \, u_{i,\pgd}^j (\bmu) + \mat{E}_{\Gamma_i \to \Omega_i} \bLambda_{\Gamma_i} \,,
\end{equation}
that performs a (much cheaper) linear combination of the functions of the PGD expansion.

Therefore, using the precomputed local surrogate models, the interface system~\eqref{eq:twoDomainInterfaceAlgebraic} can be expressed as: for any $\bar{\bmu} \in \mathcal{P}$, find $\bLambda_{\Gamma_1}$ and $\bLambda_{\Gamma_2}$ such that
\begin{equation}\label{eq:twoDomainInterfaceAlgebraicPGD}
 \Sigma^\pgd
 \begin{pmatrix}
     \bLambda_{\Gamma_1} \\[2pt]
     \bLambda_{\Gamma_2}
 \end{pmatrix} \\
 =
 \begin{pmatrix}
     \mathbf{g}_{\Gamma_1} + \mat{R}_{\Omega_2 \to \Gamma_1} u_{2,\pgd}^{0,f} (\bar{\bmu}) \\[2pt]
     \mathbf{g}_{\Gamma_2} + \mat{R}_{\Omega_1 \to \Gamma_2} u_{1,\pgd}^{0,f} (\bar{\bmu})
 \end{pmatrix},
\end{equation}
with
\begin{equation*}
\Sigma^\pgd = 
 \begin{pmatrix}
 \mat{I}_{\Gamma_1\Gamma_1} & -\mat{R}_{\Omega_2 \to \Gamma_1} \mathcal{P}_2^\pgd \\
 -\mat{R}_{\Omega_1 \to \Gamma_2} \mathcal{P}_1^\pgd & \mat{I}_{\Gamma_2\Gamma_2}
 \end{pmatrix}.
\end{equation*}

The online coupling phase thus consists of solving system~\eqref{eq:twoDomainInterfaceAlgebraicPGD} still using GMRES, but now the computational cost of each iteration is mainly due to performing the linear combinations~\eqref{eq:poissonPGD} in each subdomain $\Omega_i$ ($i=1,2$). At convergence of the GMRES iterations, we can define the global surrogate model as
\begin{equation}\label{eq:globalSoln}
	u_\pgd (\bar{\bmu}) = 
	\left\{
	\begin{array}{ll}
	 \mathcal{P}_1^\pgd \bLambda_{\Gamma_1} + u_{1,\pgd}^{0,f} (\bar{\bmu}) + g_{\Omega_1}^h (\bar{\bmu}) & \text{ in } \Omega_1, \\
	 \mathcal{P}_2^\pgd \bLambda_{\Gamma_2} + u_{2,\pgd}^{0,f} (\bar{\bmu}) + g_{\Omega_2}^h (\bar{\bmu}) & \text{ in } \Omega_2 \setminus \Omega_{12}.
	\end{array}
	\right.
\end{equation}

\section{Critical comparison of PGD-based local surrogate models}
\label{sc:comparison}

The methodology introduced in Sect.~\ref{sc:offlineOnline} is similar to the one proposed in~\cite{Discacciati:2024:CMAME}, but it presents several computational advantages.
In this section, we begin by briefly summarizing the approach introduced in~\cite{Discacciati:2024:CMAME}, and then we compare it to the method studied in this work to highlight the superior performance of the latter, both in the offline and in the online phase.

\subsection{Local surrogate models with active interface parameters}
\label{sc:pgdAIP}
\begin{figure}[!ht]
    \centering
    \resizebox{0.9\textwidth}{!}{
    \begin{tikzpicture}
        \newcommand{\intNodes}[3]{
            \draw[dashed, #3] (\longSide + #1 * \delta / 7, \ang * #1 * \delta / 7 + #2) -- (\longSide + #1 * \delta / 7, #1 * \ang * \delta / 7 );
            \fill[#3] (\longSide + #1 * \delta / 7, \ang * #1 * \delta / 7 + #2) circle (0.07 cm);
            \fill[#3] (\longSide + #1 * \delta / 7, #1 * \ang * \delta / 7 ) circle (0.07 cm);
        }

        \newcommand{\intFun}[3]{
            \draw[thick, red] (\longSide + #1 * \delta / 7, #1 * \ang * \delta / 7 + #2 ) -- (\longSide + #1 * \delta / 7 + \delta / 7, #1 * \ang * \delta / 7 +  \ang * \delta / 7 + #3);
        }

        \def\ang{1}
        \def\delta{1.5}
        \def\longSide{5}
        
        \begin{scope}
            \draw (0, 0) -- (\delta, \delta*\ang);
            \draw (0, 0) -- (\longSide, 0);
            \draw(\delta, \delta*\ang) -- (\delta + \longSide, \delta*\ang);
            \draw(\longSide, 0) -- (\delta + \longSide, \delta*\ang);
    
            \node at (\longSide / 2 + \delta / 2, \delta * \ang / 2) {$\Omega_i$};
    
            \draw[-stealth, thick] (\longSide + \delta + 0.75, \delta * \ang + 0.25) node[anchor=west] {$\Gamma_i$} parabola (\longSide + \delta + 0.05, \delta * \ang);
            
            \intNodes{1}{1.6}{blue};
            \intNodes{2}{2}{blue};
            \intNodes{3}{1.5}{blue};
    
            \intNodes{4}{1.8}{green!80!black};
            \intNodes{5}{2}{green!80!black};
            \intNodes{6}{1.4}{green!80!black};
    
            \draw [blue, decorate, decoration = {brace, mirror}, thick] (\longSide + \delta / 7 + 0.1, \ang * \delta / 7 - 0.1) -- (\longSide + 3 * \delta / 7 + 0.1, 3 * \ang * \delta / 7 - 0.1) node[pos=0.2, right=3pt, blue]{$\mathcal{N}_i^1$};
    
            \draw [green!80!black, decorate, decoration = {brace, mirror}, thick] (\longSide + 4 * \delta / 7 + 0.1, 4 * \ang * \delta / 7 - 0.1) -- (\longSide + 6 * \delta / 7 + 0.1, 6 * \ang * \delta / 7 - 0.1) node[pos=0.2, right=3pt, green!80!black]{$\mathcal{N}_i^2$};
    
            \intFun{1}{1.6}{2};
            \intFun{2}{2}{1.5};
            \intFun{3}{1.5}{1.8};
            \intFun{4}{1.8}{2};
            \intFun{5}{2}{1.4};
        \end{scope}

        \begin{scope}[shift={(-4.5, -4)}]
            \draw (0, 0) -- (\delta, \delta*\ang);
            \draw (0, 0) -- (\longSide, 0);
            \draw(\delta, \delta*\ang) -- (\delta + \longSide, \delta*\ang);
            \draw(\longSide, 0) -- (\delta + \longSide, \delta*\ang);
    
            \node at (\longSide / 2 + \delta / 2, \delta * \ang / 2) {$\Omega_i$};
    
            \draw[-stealth, thick] (\longSide + \delta + 0.75, \delta * \ang + 0.25) node[anchor=west] {$\Gamma_i$} parabola (\longSide + \delta + 0.05, \delta * \ang);
            
            \intNodes{1}{1.6}{blue};
            \intNodes{2}{2}{blue};
            \intNodes{3}{1.5}{blue};

            \intNodes{4}{0}{lightgray};
            \intNodes{5}{0}{lightgray};
            \intNodes{6}{0}{lightgray};
    
            \draw [blue, decorate, decoration = {brace, mirror}, thick] (\longSide + \delta / 7 + 0.1, \ang * \delta / 7 - 0.1) -- (\longSide + 3 * \delta / 7 + 0.1, 3 * \ang * \delta / 7 - 0.1) node[pos=0.2, right=3pt, blue]{$\mathcal{N}_i^1$};
    
            \draw [lightgray, decorate, decoration = {brace, mirror}, thick] (\longSide + 4 * \delta / 7 + 0.1, 4 * \ang * \delta / 7 - 0.1) -- (\longSide + 6 * \delta / 7 + 0.1, 6 * \ang * \delta / 7 - 0.1) node[pos=0.2, right=3pt, lightgray]{$\mathcal{N}_i^2$};
    
            \intFun{1}{1.6}{2};
            \intFun{2}{2}{1.5};
            \intFun{3}{1.5}{0};
            \intFun{4}{0}{0};
            \intFun{5}{0}{0};
        \end{scope}

        \begin{scope}[shift={(3, -4)}]
            \draw (0, 0) -- (\delta, \delta*\ang);
            \draw (0, 0) -- (\longSide, 0);
            \draw(\delta, \delta*\ang) -- (\delta + \longSide, \delta*\ang);
            \draw(\longSide, 0) -- (\delta + \longSide, \delta*\ang);
    
            \node at (\longSide / 2 + \delta / 2, \delta * \ang / 2) {$\Omega_i$};
    
            \draw[-stealth, thick] (\longSide + \delta + 0.75, \delta * \ang + 0.25) node[anchor=west] {$\Gamma_i$} parabola (\longSide + \delta + 0.05, \delta * \ang);
            
            \intNodes{1}{0}{lightgray};
            \intNodes{2}{0}{lightgray};
            \intNodes{3}{0}{lightgray};

            \intNodes{4}{1.8}{green!80!black};
            \intNodes{5}{2}{green!80!black};
            \intNodes{6}{1.4}{green!80!black};
    
            \draw [lightgray, decorate, decoration = {brace, mirror}, thick] (\longSide + \delta / 7 + 0.1, \ang * \delta / 7 - 0.1) -- (\longSide + 3 * \delta / 7 + 0.1, 3 * \ang * \delta / 7 - 0.1) node[pos=0.2, right=3pt, lightgray]{$\mathcal{N}_i^1$};
    
            \draw [green!80!black, decorate, decoration = {brace, mirror}, thick] (\longSide + 4 * \delta / 7 + 0.1, 4 * \ang * \delta / 7 - 0.1) -- (\longSide + 6 * \delta / 7 + 0.1, 6 * \ang * \delta / 7 - 0.1) node[pos=0.2, right=3pt, green!80!black]{$\mathcal{N}_i^2$};
    
            \intFun{1}{0}{0};
            \intFun{2}{0}{0};
            \intFun{3}{0}{1.8};
            \intFun{4}{1.8}{2};
            \intFun{5}{2}{1.4};
        \end{scope}
        
        \draw[- stealth] (1.9375, -0.3) -- (-0.1875 , -2);
        \draw[- stealth] (4, -0.3) -- (5.3 , -2);
        
        \node at ( 5.2,  3.0) {$\lambda_i(\bx;\bmu)$};
        \node at ( 2.8, -1.8) {$\sum_{q \in \mathcal{N}_i^1} \Lambda_i^q(\bmu)\, \eta_i^q (\bmu)$};
        \node at ( 7.4, -1.0) {$\sum_{q \in \mathcal{N}_i^2} \Lambda_i^q(\bmu)\, \eta_i^q (\bmu)$};

    \end{tikzpicture}
  }
    \caption{Example of clustering of the interface nodes into two sets of active interface parameters $\mathcal{N}_i^1$ and $\mathcal{N}_i^2$.}
    \label{fig:lambdaSplittingActive}
\end{figure}
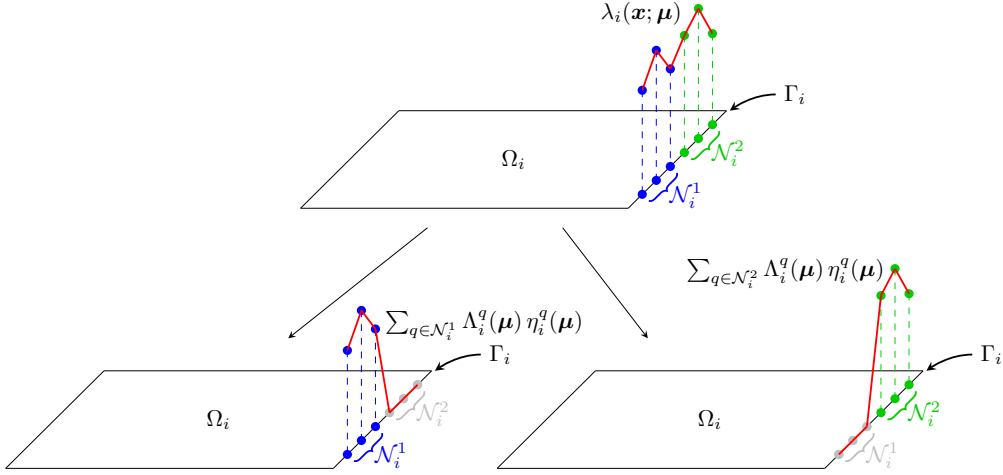

In \cite{Discacciati:2024:CMAME}, the same representations~\eqref{eq:lambdaRepresentation} and~\eqref{eq:lambdaExtension} of the arbitrary interface function $\lambda_i^h(\bmu)$ and of its extension $\lambda_{\Omega_i}^h(\bmu)$ are considered. However, to construct the local functions $u_{i,h}^{0,\lambda}(\bmu)$, problem~\eqref{eq:galerkinLambda} is addressed directly, using $\lambda_{\Omega_i}^h(\bmu)$ at the right-hand side and handling the unknown coefficients $\{ \Lambda_i^j(\bmu) \}_{j=1,\ldots,N_{\Gamma_i}}$ as additional auxiliary problem parameters.
This increases the dimensionality of the local parametric problems~\eqref{eq:galerkinLambda}, and can become problematic if $N_{\Gamma_i} \gg 1$. In fact, it is well known that PGD, like any other model order reduction method, cannot efficiently handle large numbers of independent parameters. Therefore, to overcome the difficulty of handling too many interface parameters $\Lambda_i^j(\bmu)$, these are clustered in $N_i$ sufficiently small, disjoint sets $\mathcal{N}_i^j$ ($j=1,\ldots,N_i$) of so-called \textit{active interface parameters} with $\bigcup_{j=1,\ldots,N_i} \mathcal{N}_i^j = \{1,\ldots,N_{\Gamma_i}\}$ and with $\text{card}(\mathcal{N}_i^j) \ll N_{\Gamma_i}$. Then,~\eqref{eq:lambdaRepresentation} is replaced by 
\begin{equation}\label{eq:splittingBoundaryParameters}
\lambda_i^h (\bx; \bmu) = 
\sum_{q \in \mathcal{N}^1_i} \Lambda^q_i (\bmu) \, \eta^q_i (\bx) + 
\sum_{q \in \mathcal{N}^2_i} \Lambda^q_i (\bmu) \, \eta^q_i (\bx) + \ldots +
\sum_{q \in \mathcal{N}^{N_i}_i} \Lambda^q_i (\bmu) \, \eta^q_i (\bx) .
\end{equation}
A schematic representation of the splitting~\eqref{eq:splittingBoundaryParameters} is shown in Fig.~\ref{fig:lambdaSplittingActive}, which highlights the clustering of the nodes in sets of active interface parameters, differently from the novel approach proposed in this work and illustrated in Fig.~\ref{fig:lambdaSplitting}. 

In the offline phase, a surrogate model depending on the problem parameters $\bmu$ and on the newly introduced active interface parameters $\{ \Lambda^q_i (\bmu) \}_{q \in \mathcal{N}_i^j}$ is thus constructed.
The $N_{\Gamma_i}$ local independent problems~\eqref{eq:unitaryProblems} are replaced by the following $N_i < N_{\Gamma_i}$ local independent problems: for $j=1,\ldots,N_i$, for all $\bmu \in \mathcal{P}$, and for all $\bLambda_i^j = \{\Lambda_i^q(\bmu) \}_{q\in \mathcal{N}_i^j} \in \mathcal{Q}_i^j$, find $\widehat{u}_{i,h}^j (\bmu, \bLambda_i^j) \in V_i^h$ such that
\begin{equation}\label{eq:boundaryProbsCluster}
\mathcal{A}_i (\widehat{u}_{i,h}^j (\bmu,\bLambda_i^j), v_{i,h}; \bmu) = - \mathcal{A}_i \left( \sum_{q \in \mathcal{N}_i^j} \Lambda_i^q(\bmu) \eta^q_i (\bx), v_{i,h}; \bmu\right) \qquad \forall v_{i,h} \in V_i^h ,
\end{equation}
where $\mathcal{Q}_i^j = \bigtimes_{q \in \mathcal{N}_i^j} \mathcal{J}_i^q$ is the space of the auxiliary parametric interface parameters, with each $\mathcal{J}_i^q \subset \mathbb{R}$ ($q = 1, \dots , N_{\Gamma_i}$) being a compact set.

Exploiting the linearity of the parametric PDE under analysis, but using now the splitting~\eqref{eq:splittingBoundaryParameters} (instead of~\eqref{eq:lambdaRepresentation}), one finally obtains the following representation of $u_{i,h}^{0,\lambda}$ (instead of~\eqref{eq:ulambdaExpansion}):
\begin{equation}\label{eq:solutionSplitCluster}
u_{i,h}^{0,\lambda}(\bmu) = \sum_{j=1}^{N_i} \widehat{u}_{i,h}^j(\bmu, \bLambda_i^j) .
\end{equation}
The corresponding PGD approximation of $\widehat{u}_{i,h}^j(\bmu, \bLambda_i^j)$ is sought in the form 
\begin{equation}\label{eq:clusterBasisModes}
\widehat{u}_{i,h}^j(\bmu, \bLambda_i^j) \approx
\widehat{u}_{i,\pgd}^j (\bmu,\bLambda_i^j) = \sum_{m=1}^{\widehat{M}_i^j} \widehat{U}_{i,j}^m(\bx)\, \widehat{\phi}_{i,j}^m (\bmu) \, \psi_{i,j}^m (\bLambda_i^j),
\end{equation}
where $\widehat{U}_{i,j}^m(\bx)$ and $\widehat{\phi}_{i,j}^m (\bmu)$ still represent the $m$th spatial and parametric modes, respectively, while $\psi_{i,j}^m (\bLambda_i^j)$ is an additional mode for the auxiliary interface parameters $\bLambda_i^j$. While the spatial modes are discretized by finite elements, both parametric modes are discretized by pointwise collocation.

Given the separated form~\eqref{eq:PGDoperator} of the PGD operators, the local surrogate models $\widehat{u}_{i,\pgd}^j (\bmu,\bLambda_i^j)$, for $j = 1, \dots, N_i$, are then computed by solving the local problems
\begin{equation}\label{eq:problemLocalPGDClustering}
 \mathcal{A}_i^{\pgd}(\widehat{u}_{i,\pgd}^j (\bmu,\bLambda_i^j), v_{i,h}; \bmu) = - \mathcal{A}_i^{\pgd} \left(\sum_{q \in \mathcal{N}_i^j} \Lambda_i^q(\bmu) \eta^q_i (\bx), v_{i,h} ; \bmu \right)\,,
\end{equation}
for all $v_{i,h} \in V_i^h$, $\bmu \in \mathcal{P}$, and $\bLambda_i^j = \{ \Lambda_i^q(\bmu) \}_{q\in \mathcal{N}_i^j} \in \mathcal{Q}_i^j$. 

To retrieve the parametric solution~\eqref{eq:solSuperposition} of problem~\eqref{eq:galerkinTwoDomain}, the functions $\{\widehat{u}_{i,\pgd}^j (\bmu,\bLambda_i^j)\}_{j = 1, \dots, N_i}$ obtained by solving~\eqref{eq:problemLocalPGDClustering} are combined with the function ${u}^{0,f}_{i,\pgd}(\bmu)$ whose computation is described in Sect.~\ref{sc:offlinePhase} and is the same as in~\cite{Discacciati:2024:CMAME}.
This provides the PGD-based local surrogate model for the parametric problem~\eqref{eq:globalProblem} in subdomain $\Omega_i$, using the clustering approach for active interface parameters.

Finally, the online phase in~\cite{Discacciati:2024:CMAME} follows the same ideas as in Sect.~\ref{sc:onlinePhase}, but the local PGD operator $\mathcal{P}_i^\pgd$ in~\eqref{eq:poissonPGD} is replaced by
\begin{equation}\label{eq:poissonPGDclustered}
\widehat{\mathcal{P}}_i^\pgd : \bLambda_{\Gamma_i} \to \sum_{j=1}^{N_i} \widehat{u}_{i,\pgd}^j (\bmu, \bLambda_i^j) + \mat{E}_{\Gamma_i \to \Omega_i} \bLambda_{\Gamma_i}\,,
\end{equation}
so that, at each GMRES iteration, one must extend the nodal values $\bLambda_{\Gamma_i}$ in $\Omega_i$, evaluate $\widehat{u}_{i,\pgd}^j (\bmu, \bLambda_i^j)$ at the given values of $\bLambda_{\Gamma_i}$, and perform the sum of the functions of the PGD expansion.

\subsection{Complexity assessment of offline and online phase}
\label{sc:Assessment}

This section presents a critical discussion of the advantages and disadvantages of the proposed local surrogate model with reduced dimensionality presented in Sect.~\ref{sc:OfflineRedDim}, showcasing its superior performance, both in the offline and in the online phase, with respect to the active interface parameters strategy introduced in~\cite{Discacciati:2024:CMAME}.

Although in the offline phase the approach in~\cite{Discacciati:2024:CMAME} requires solving a smaller number $N_i < N_{\Gamma_i}$ of local problems~\eqref{eq:problemLocalPGDClustering}, each computation depends not only on the problem parameters $\bmu$, but also on the auxiliary interface parameters $\bLambda_i^j$.
On the contrary, problems~\eqref{eq:problemLocalPGDUnity} only involve the parameter $\bmu$, thus retaining the same dimensionality as the original equation~\eqref{eq:globalProblem}.
Hence, despite the larger number of local surrogate models to be computed, the fact that all problems~\eqref{eq:problemLocalPGDUnity} are independent from one another, that they can be easily solved in parallel, and that they have a reduced dimensionality, makes the computational cost of obtaining the PGD expansion much cheaper than in the case of clustered nodes presented in~\cite{Discacciati:2024:CMAME}.

Moreover, before solving~\eqref{eq:problemLocalPGDClustering}, a suitable interface parametric space $\mathcal{Q}_i^j$ must be identified to represent an arbitrary interface function that is non-zero at the nodes whose indices belong to the set of active interface parameters $\mathcal{N}_i^j$. In principle, each interval $\mathcal{J}_i^q$ depends on $\bmu$ and it should be defined considering physical information about the solution of problem~\eqref{eq:globalProblem} (e.g., the maximum and minimum values that $u(\bmu)$ can attain within $\Omega$). In the approach proposed in this work, this is not needed, as the auxiliary problems~\eqref{eq:problemLocalPGDUnity} are independent of $\bLambda_i^j$.
It should also be noted that not all possible combinations of the parameters $\{ \Lambda_i^q(\bmu) \}_{q \in \mathcal{N}_i^j}$ considered in~\eqref{eq:problemLocalPGDClustering} may actually be significant for the problem at hand. For instance, taking very different values of $\Lambda_i^q$ at adjacent nodes on $\Gamma_i$ could represent a highly oscillatory function on $\Gamma_i$, which is unlikely to be relevant for the solution of an elliptic problem in a smooth domain with sufficiently regular data. Solving boundary value problems with such localized oscillatory features can become challenging for PGD due to the possibly large number of modes needed to correctly represent the behavior of the solution. Thus, obtaining all such modes can unnecessarily increase the computational cost of the offline phase. 

Finally, there is no need to cluster the interface nodes when working with~\eqref{eq:problemLocalPGDUnity}, so that the practical computer implementation of the approach presented in this work is much simpler.

\smallskip

Concerning the online phase, the methods presented in this work and in~\cite{Discacciati:2024:CMAME} differ in the definition of the local PGD operator, see~\eqref{eq:poissonPGD} for the interface problem with reduced dimensionality and~\eqref{eq:poissonPGDclustered} for the case of active interface parameters.
Given a value, say, $\bar{\bmu}$, of the problem parameter, the former approach first evaluates all local surrogate models at $\bar{\bmu}$. Note that this is done before starting the GMRES iterations to solve~\eqref{eq:twoDomainInterfaceAlgebraicPGD}, as a preprocess for the online evaluation. Then, at the $k$th GMRES iteration, given any array of values $\bLambda_{\Gamma_i}^{(k)}$, \eqref{eq:poissonPGD} requires to perform the linear combination
\begin{equation}\label{eq:linCombOnline}
\sum_{j=1}^{N_{\Gamma_i}} (\Lambda_i^j(\bar{\bmu}))^{(k)} \, u_{i,\pgd}^j(\bar{\bmu})
\end{equation}
and the extension $\mat{E}_{\Gamma_i \to \Omega_i} \bLambda_{\Gamma_i}^{(k)}$. Thus, the computational cost of each GMRES iteration is associated with these operations: for each subdomain, the linear combination~\eqref{eq:linCombOnline} involving $N_{\Gamma_i}$ terms and the extension.
The pre-evaluation of the PGD expansion at $\bar{\bmu}$ before the GMRES iterations is not possible in the clustered approach of~\cite{Discacciati:2024:CMAME} with the active interface parameters. Indeed, recalling~\eqref{eq:clusterBasisModes}, \eqref{eq:poissonPGDclustered} requires the evaluation of
\begin{equation}\label{eq:poissonPGDclusteredEvaluation}
\sum_{j=1}^{N_i} \widehat{u}_{i,\pgd}^j (\bar{\bmu}, (\bLambda_i^j)^{(k)}) = \sum_{j = 1}^{N_i} \sum_{m = 1}^{\widehat{M}_i^j} \widehat{U}_{i, j}^m(\bx) \, \widehat{\phi}_{i, j}^m(\bar{\bmu}) \, \psi_{i, j}^m((\bLambda_i^j)^{(k)})
\end{equation}
and the extension $\mat{E}_{\Gamma_i \to \Omega_i} \bLambda_{\Gamma_i}^{(k)}$.
Note that evaluating $\widehat{u}_{i,\pgd}^j (\bar{\bmu}, (\bLambda_i^j)^{(k)})$ requires $(\bLambda_{i}^j)^{(k)}$ to be known at all DOFs of the cluster $\mathcal{N}_i^j$.
Hence, at each GMRES iteration, the double summation on the right-hand side of~\eqref{eq:poissonPGDclusteredEvaluation} must be performed in each subdomain, significantly increasing the computational cost compared to the case described in~\eqref{eq:linCombOnline}.

Moreover, it must be noticed that the values $\bLambda_{\Gamma_i}^{(k)}$ at the DOFs of the interface $\Gamma_i$ are generated by GMRES without any control from the user. This may imply that at least some of the coefficients of $\bLambda_{\Gamma_i}^{(k)}$ may not coincide with those obtained when discretizing the parametric space $\mathcal{Q}_i^j$ by collocation, and they may even fall outside of the parametric space $\mathcal{Q}_i^j$. In the latter case, the convergence of the GMRES iterations may be completely jeopardized, while in the former a linear interpolation of the parametric modes $\psi_{i,j}^m$ associated with the available discrete values in $\mathcal{Q}_i^j$ closest to $\bLambda_{\Gamma_i}^{(k)}$ must be performed. Although this operation is not computationally demanding, its repeated execution when applying the PGD operator $\widehat{\mathcal{P}}_i^\pgd$ leads to a higher computational cost than the one required by operator $\mathcal{P}_i^\pgd$.
On the contrary, this cannot occur using the interface problem with reduced dimensionality proposed in Sect.~\ref{sc:offlineOnline}. Indeed, it is always possible to compute the linear combination~\eqref{eq:linCombOnline} for any given coefficients $\bLambda_{\Gamma_i}^{(k)}$ once the PGD expansions $\{u^{j}_{i,\pgd}(\bmu)\}_{j = 1, \dots, N_{\Gamma_i}}$ are available. The resulting online coupling algorithm is thus more robust and computationally less expensive than the one introduced in~\cite{Discacciati:2024:CMAME}, as showcased by the numerical experiments in Sect.~\ref{sc:numericalResults}.

\section{Numerical results}
\label{sc:numericalResults}
In this section, we provide numerical experiments\footnote{The numerical results presented in this section have been obtained using a PC with CPU Intel$^{\mbox{\tiny{\textregistered}}}$ Core\texttrademark\; i5-11400 @ 2.60GHz and 8GB RAM.} benchmarking the local surrogate model with reduced dimensionality presented in Sect.~\ref{sc:offlineOnline} with the strategy featuring active interface parameters introduced in~\cite{Discacciati:2024:CMAME} and recalled in Sect.~\ref{sc:pgdAIP}. 
For all test cases, in the offline phase, the local parametric problems~\eqref{eq:problemLocalPGDData},~\eqref{eq:problemLocalPGDUnity}, and~\eqref{eq:problemLocalPGDClustering} are solved using the Encapsulated PGD Algebraic Toolbox~\cite{Diez:2020:ACME} with tolerance $10^{-4}$ to stop the PGD enrichment process, and a compression algorithm~\cite{DM-MZH:15} with tolerance $10^{-3}$ is applied to eliminate redundant modes. In the online phase, the interface system~\eqref{eq:twoDomainInterfaceAlgebraicPGD} is solved by GMRES with stopping tolerance $10^{-6}$ on the relative residual.

\subsection{Sensitivity to the dimensionality of the local surrogate models}
\label{sc:sensitivityTesting}
In this section, the effect of the dimensionality of the local surrogate model is studied by considering different descriptions of the trace of the solution at the interface, namely,~\eqref{eq:lambdaRepresentation} and~\eqref{eq:splittingBoundaryParameters}, the latter varying the number of active interface parameters from 1 to 5.

To perform the sensitivity study, a two-dimensional synthetic test case with analytical solution and $N_p=1$ is considered.
Problem~\eqref{eq:globalProblem} is set in the domain $\Omega = (0,2)\times (0,1)$, split into two overlapping subdomains $\Omega_1 = (0, 1.05) \times (0, 1)$ and $\Omega_2 = (0.95, 2) \times (0, 1)$.
The parametric diffusion coefficient is defined as $\nu(\mu) = 1+\mu x$, with $\mu \in [1,50]$ being a scalar parameter, and the source term $f(\mu)$ is given by
\begin{equation*}
\begin{aligned}
f(\mu) =& 8\pi^2 \sin(2\pi x) \, \sin(2\pi y) \\
&+ \mu [2\pi(4\pi x \sin(2\pi x) - \cos(2\pi x))\sin(2\pi y) - x(x-2) - y(y-1)] \\
&+ \mu^2 [y(y-1)(1-2x) - x^2(x-2)] .
\end{aligned}
\end{equation*}
The corresponding exact solution is
\begin{equation*}
u_{\text{ex}}(\mu) = \sin(2\pi x) \, \sin(2\pi y) + \frac{\mu}{2} \, xy(y-1)(x-2) ,
\end{equation*}
and homogeneous Dirichlet boundary conditions are applied on the entire boundary $\partial\Omega$.

A uniform structured mesh of quadrilateral elements is defined in the domain $\Omega$, with mesh size $h=5 \times 10^{-2}$, and the local meshes of the two subdomains coincide in the overlap of width $2h$. 
The parametric domain is discretized using a uniform mesh of size $h_\mu = 10^{-3}$.
A continuous Galerkin finite element approximation with $\mathbb{Q}_1$ Lagrange basis functions is used for the spatial discretization, whereas pointwise collocation is employed in the parametric direction.

Given this discretization in space, there are 19 nodes on each interface.
Four configurations of PGD-based surrogate models are analyzed in this section, differing in the strategy employed to handle the large number of parameters introduced by the definition of the auxiliary trace variable.
On the one hand, the strategy proposed in Sect.~\ref{sc:offlineOnline} imposes unitary interface conditions independently at each node on the interface, leading to a number of low-dimensional subproblems equal to the number of DOFs on the interface, namely 19.
On the other hand, clusters of $\nAIP$ active interface parameters (namely, 1, 3, and 5) are considered for the approach described in Sect.~\ref{sc:pgdAIP}, lowering the number of local subproblems while increasing their dimensionality. 
Moreover, all configurations require the solution of $\nDP=1$ additional problem~\eqref{eq:problemLocalPGDData} with parametric data in each subdomain.

Let $\nIP$ denote the number of interface problems~\eqref{eq:problemLocalPGDUnity} or~\eqref{eq:problemLocalPGDClustering} and $\dIP$ the number of dimensions of each subproblem defined as $\dIP=d+N_p+\nAIP$.
Table~\ref{tab:sensitivityTestingOffline} presents the details of the offline phase for each configuration of the local surrogate model.
Although the strategy based on unitary interface conditions (Sect.~\ref{sc:offlineOnline}) requires solving the largest number of interface problems, the introduction of active interface parameters (Sect.~\ref{sc:pgdAIP}) results in a significant increase in the complexity of the algorithm. Indeed, by simply parametrizing the nodal value of the interface condition (that is, increasing $\dIP$ from 3 to 4), the corresponding offline computing time $\Toff$ increases from $5.66$ to $9.80$~s.
This growth becomes even more evident for $\nAIP=3$ (respectively, 5), with the corresponding interface problems being of dimension 6 (respectively, 8) and the offline execution requiring $112.14$~s (respectively, $482.68$~s).
Hence, the PGD algorithm to construct the local surrogate models tends to suffer as the number of involved dimensions $\dIP$ increases and the strategy proposed in Sect.~\ref{sc:offlineOnline} thus outperforms the active interface parameters approach. This is also testified by the growth in the number of computed modes for $\nAIP=5$, showcasing that redundant information is determined by the PGD in this case, likely associated with inadmissible combinations of parameters (e.g., featuring large node-to-node variations of the solution).
\begin{table}[!htb]
\centering
\begin{tabular}{| l | c || c | c | c || c | c || c |}
\hline
Surrogate & \multirow{2}{*}{$\nAIP$} & \multirow{2}{*}{$\nDP$} & \multirow{2}{*}{$\nIP$} & \multirow{2}{*}{$\dIP$} & \multicolumn{2}{c||}{Modes$^\dagger$} &\multirow{2}{*}{$\Toff$ (s)} \\
\cline{6-7}
strategy & & & & & $\Omega_1$ & $\Omega_2$ & \\
\hline
Sect.~\ref{sc:offlineOnline} & - & 1 & 19 & 3 & 68 (106) & 56 (62) & 5.66 \\
\hline
\multirow{3}{*}{Sect.~\ref{sc:pgdAIP}} & 1 & 1 & 19 & 4 & 68 (106) & 56 (62) & 9.80 \\
                        & 3 & 1 & 7 & 6 & 59 (114) & 40 (62) & 112.14 \\
                        & 5 & 1 & 4 & 8 & 56 (137) & 34 (73) & 482.68 \\
\hline
\end{tabular}

\caption{Offline phase for the local surrogate models using four different approaches to handle the interface parameters. $^\dagger$The number in brackets denotes the number of modes before compression.}
\label{tab:sensitivityTestingOffline}
\end{table}

The local surrogate model with reduced dimensionality presents computational advantages also in the online phase, as showcased in Table~\ref{tab:sensitivityTestingOnline} for two values of the parameter $\mu$.
The model of Sect.~\ref{sc:offlineOnline} displays faster convergence, reducing the number of GMRES iterations to 9 from the 18 required to achieve convergence in the case of 5 active interface parameters.
It is worth noticing that the corresponding online computing times $\Ton$ are reduced by more than a factor 10. Indeed, this computational gain is the result of both reducing the number of GMRES iterations and reducing the cost of the evaluation of the surrogate model, avoiding the interpolation of the nodal values of the interface condition (see equations~\eqref{eq:linCombOnline} and~\eqref{eq:poissonPGDclusteredEvaluation}).
\begin{table}[!htb]
\centering
\begin{tabular}{| l | c || c | c || c | c |}
\hline
Surrogate & \multirow{2}{*}{$\nAIP$} & \multicolumn{2}{c||}{$\nGMRES$} & \multicolumn{2}{c|}{$\Ton$ (s)} \\
\cline{3-6}
strategy & & $\mu=3$ & $\mu=30$ & $\mu=3$ & $\mu=30$ \\
\hline
Sect.~\ref{sc:offlineOnline} & - & 9 & 9 & 0.11 & 0.11 \\
\hline
\multirow{3}{*}{Sect.~\ref{sc:pgdAIP}} & 1 & 13 & 11 & 0.94 & 0.79 \\
                        & 3 & 16 & 15 & 1.24 & 1.17 \\
                        & 5 & 18 & 14 & 1.49 & 1.24 \\
\hline
\end{tabular}

\caption{Online phase for the local surrogate models using four different approaches to handle the interface parameters.}
\label{tab:sensitivityTestingOnline}
\end{table}

\begin{figure}[!htb]
    \centering
    \subfigure[Surrogate with reduced dimensionality]{\includegraphics[width=0.4\columnwidth]{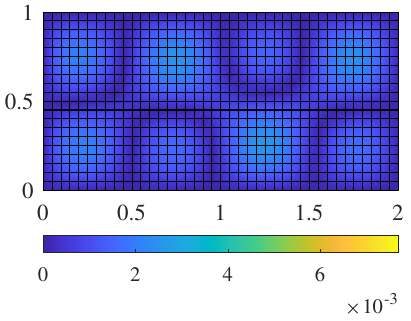}}
    \hspace{10pt}
    \subfigure[Surrogate with $\nAIP=1$]{\includegraphics[width=0.4\columnwidth]{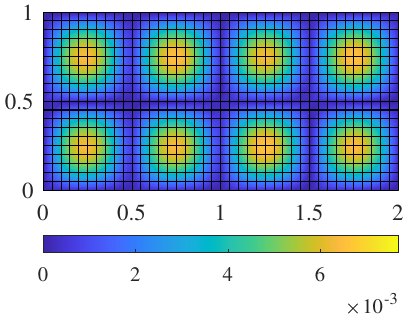}}
    
    \subfigure[Surrogate with $\nAIP=3$]{\includegraphics[width=0.4\columnwidth]{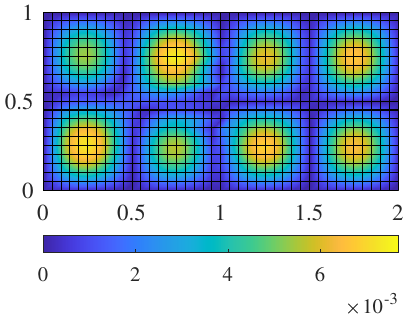}}
    \hspace{10pt}
    \subfigure[Surrogate with $\nAIP=5$]{\includegraphics[width=0.4\columnwidth]{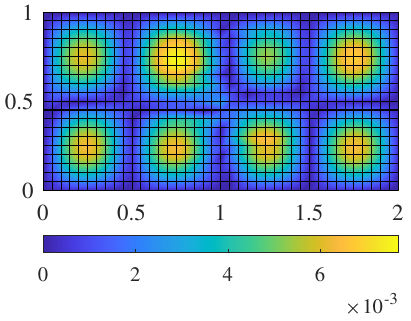}}
    
    \caption{Map of the scaled nodal error $| u_{\pgd}(\mu) - u_{\text{ex}}(\mu)|/\max_{\Omega} |u_{\text{ex}}(\mu)|$ for $\mu=3$ using four different approaches to handle the interface parameters.}
    \label{fig:sensitivityTestingErrors}
\end{figure}

Let $E^{\pgd}_{2}$ be the relative error, measured in the $L^2(\Omega)$ norm, of the PGD surrogate model $u_{\pgd}(\mu)$ with respect to the exact solution $u_{\text{ex}}(\mu)$ for a fixed value of $\mu$, that is, 
\begin{equation}\label{eq:errL2Def}
E^{\pgd}_{2} = \frac{\| u_{\pgd}(\mu) - u_{\text{ex}}(\mu) \|_{L^2(\Omega)}}{\| u_{\text{ex}}(\mu) \|_{L^2(\Omega)}} .
\end{equation}
To assess the accuracy of the reconstructed global solutions using the four surrogate models, we compute $E^{\pgd}_{2}$ for selected values of the parameter $\mu$.
The four configurations attain comparable relative errors of $9 \times 10^{-3}$ for $\mu=3$ and $3 \times 10^{-3}$ for $\mu = 30$, achieving the same accuracy as a reference full-order solution computed using the finite element method on the entire domain.
Finally, Fig.~\ref{fig:sensitivityTestingErrors} reports the map of the scaled nodal error $| u_{\pgd}(\mu) - u_{\text{ex}}(\mu)|/\max_{\Omega} |u_{\text{ex}}(\mu)|$. It should be noted that, although of the same magnitude, the error of the surrogate model with reduced dimensionality is slightly lower than the one achieved by the strategy based on active interface parameters, thus outperforming the latter approach both in terms of efficiency and accuracy.

To summarize, the proposed approach of local surrogate models with reduced dimensionality is more efficient than the original strategy introduced in~\cite{Discacciati:2024:CMAME}, independently of the number of clustered nodes selected in the latter as the set of active interface parameters.
Offline computing times are reduced of approximately 20 times with respect to the case of 3 active interface parameters and up to 85 times for $\nAIP=5$, while slightly improving the accuracy both in the overlap region and in the entire domain.
Moreover, the online phase requires fewer GMRES iterations and achieves an average speed-up between 8 and 12 times, both advantages being attributable to the absence of any interpolation procedure for the values of the interface parameters.
In the rest of this article, the surrogate model with reduced dimensionality will only be compared to the case of active interface parameters originally presented in~\cite{Discacciati:2024:CMAME}.

\subsection{Computational comparison of the offline phase}
\label{sc:compareOffline}
To further evaluate the computational gains provided by the local surrogate models with reduced dimensionality during the offline phase, we consider a more challenging problem with two parameters $\bmu=(\mu_1,\mu_2)^T$, one controlling material properties and one controlling geometry. Specifically, the convection-diffusion equation for the Poiseuille-Graetz flow in a geometrically parametrized domain, see~\cite{Pacciarini:2014:CMAME}, is studied.

Consider the domain $\Omega(\bmu) = (0, 1 + \mu_2) \times (0, 1)$ describing a channel with boundary $\partial\Omega(\bmu) = \Gamma^{D,1}(\bmu) \cup \Gamma^{D,2}(\bmu) \cup \Gamma^{N}(\bmu)$ and the different portions being disjoint by pairs.
Let $\nu(\bmu) = \mu_1^{-1}$ be the diffusion coefficient and $\balpha = (\alpha_1,\alpha_2)^T = (4y(1 - y) , 0 )^T$ the convective velocity field.
Given the parametric domain $\mathcal{P} = \mathcal{I}_1 \times \mathcal{I}_2$ for $\bmu$, with $\mathcal{I}_1 = [10^4, 2 \times 10^4]$ and $\mathcal{I}_2 = [0.5, 4]$, the convection-diffusion equation under analysis models the evolution of the temperature field inside the channel $\Omega(\bmu)$, knowing that the walls $\Gamma^{D, 1}(\bmu)$ and $\Gamma^{D, 2}(\bmu)$ are maintained at different temperatures and $\Gamma^N(\bmu)$ is an adiabatic surface. More precisely, the problem is
\begin{equation}\label{eq:paramConvDiff}
    \begin{array}{rcll}
    -\mu_1^{-1} \,\Delta u(\bmu) + \balpha \cdot \nabla u(\bmu) &=& 0 &\quad \text{in } \Omega(\bmu)\, ,\\
         u(\bmu) &=& 0 &\quad \text{on } \Gamma^{D, 1}(\bmu)\, ,\\
         u(\bmu) &=& 1 &\quad \text{on } \Gamma^{D, 2}(\bmu)\, ,\\
         \mu_1^{-1}\nabla u(\bmu)\cdot \bn(\bmu) &=& 0 &\quad \text{on } \Gamma^N(\bmu) \, ,
    \end{array}
\end{equation}
where the boundaries are
\begin{align*}
    \Gamma^{D, 1}(\bmu) &= \Gamma^{D, 1} = \{0\} \times [0, 1] \cup [0, 1] \times \{0\} \cup [0, 1] \times \{1\} \,,\\
    \Gamma^{D, 2}(\bmu) &= [1, 1 + \mu_2] \times \{0\} \cup [1, 1 + \mu_2] \times \{1\}\, ,\\
    \Gamma^N(\bmu) &= \{1 + \mu_2\} \times [0, 1] \, .
\end{align*}

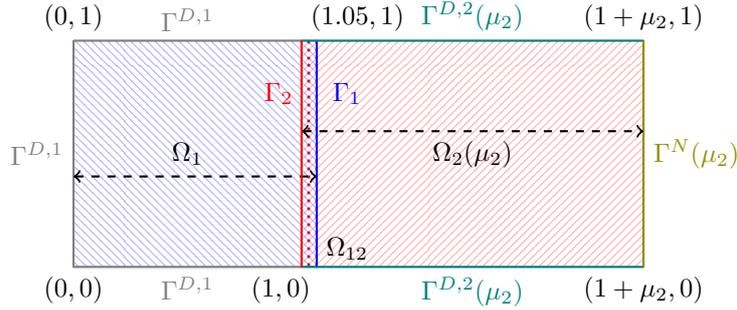
\begin{figure}[!htb]
    \centering
    \begin{tikzpicture}
        \draw[pattern=north west lines, pattern color=blue!30] (0, 0) rectangle (3.2, 3);
        \draw[pattern=north east lines, pattern color=red!30] (3, 0) rectangle (7.5, 3);        
		\draw[pattern=dots, pattern color=purple] (3, 0) rectangle (3.2, 3);
        \node at (3.6, 0.25) {$\Omega_{12}$};

        \node[red] at (2.7, 2.3) {$\Gamma_2$};
        \node[blue] at (3.6, 2.3) {$\Gamma_1$};
        
        \draw[gray, thick] (0, 0) node[black, anchor = north]{$(0,0)$} -- (1.5, 0) node[anchor = north]{$\Gamma^{D, 1}$} -- (3, 0) node[black, anchor = north east, shift={(0.25, 0)}]{$(1, 0)$} ;
        \draw[gray, thick] (0, 0) -- (0, 1.5) node[anchor = east]{$\Gamma^{D, 1}$} -- (0, 3)  node[black, anchor = south]{$(0, 1)$};
        \draw[gray, thick] (0, 3) -- (1.5, 3) node[anchor = south]{$\Gamma^{D, 1}$} -- (3, 3);
        \draw[<->, thick, dashed] (0, 1.2) -- (3.2, 1.2);
        \draw[blue, thick] (3.2, 0) -- (3.2, 3);
        \node[black, anchor = south west] at (3,3) {$(1.05, 1)$};
        
        \node at (1.5, 1.5) {$\Omega_1$};

        \draw[teal, thick] (3, 0) -- (5.25, 0) node[anchor = north]{$\Gamma^{D, 2}(\mu_2)$} --(7.5, 0) node[black, anchor = north]{$(1 + \mu_2, 0)$};
        \draw[teal, thick] (3, 3) -- (5.25, 3) node[anchor = south]{$\Gamma^{D, 2}(\mu_2)$} --(7.5, 3) node[black, anchor = south]{$(1 + \mu_2, 1)$};
        \draw[olive, thick] (7.5, 0) -- (7.5, 1.5) node[anchor = west]{$\Gamma^N(\mu_2)$} -- (7.5, 3);
        \draw[<->, thick, dashed] (3, 1.8) -- (7.5, 1.8);
        \draw[red, thick] (3, 0) -- (3, 3); 

        \node at (5.25, 1.5) {$\Omega_2(\mu_2)$};
    \end{tikzpicture}
    \caption{Computational domain for the geometrically-parametrized convection-diffusion problem.}
    \label{fig:PGflowDom}
\end{figure}

Two overlapping subdomains are introduced: a fixed subdomain $\Omega_1 = [0, 1.05] \times [0, 1]$ and a parametric subdomain $\Omega_2(\bmu) = [1, 1 + \mu_2] \times [0, 1]$, whose length is controlled by the parameter $\mu_2$ (see Fig.~\ref{fig:PGflowDom}). 
Following~\cite{Ammar-AHCCL-14}, problem~\eqref{eq:paramConvDiff} is rewritten using a reference domain.
To this end, the mapping $\Map : \hOmega_2 \times \mathcal{I}_2 \rightarrow \Omega_2(\bmu)$ is introduced to transform the coordinates $(\hx,\hy)$ of the reference domain $\hOmega_2 = [0,1] \times [0,1]$ into the physical coordinates $(x,y)$ as 
\begin{equation}\label{eq:transformation}
x = 
\begin{cases}
    1+ \hx \quad & \text{for } \hx \leq \bar{h}\,,\\[0.5em]
    \displaystyle\frac{1 - \bar{h} \hx}{1-\bar{h}} + \mu_2 \frac{\hx - \bar{h}}{1-\bar{h}}\quad &\text{for } \hx > \bar{h} \, ,
\end{cases}
\qquad
y = \hy \, ,
\end{equation}
with $\bar{h} = 5 \times 10^{-2}$.

\begin{figure}[!htb]
    \centering
    \subfigure[$\Omega_1$]{\includegraphics[width=0.35\columnwidth]{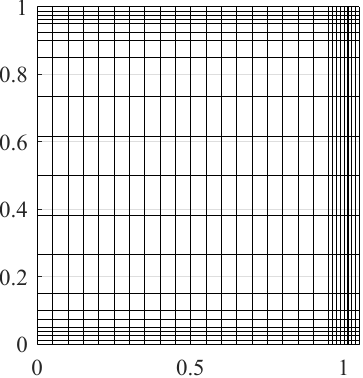}}
    \hspace{10pt}
    \subfigure[$\hOmega_2$]{\includegraphics[width=0.35\columnwidth]{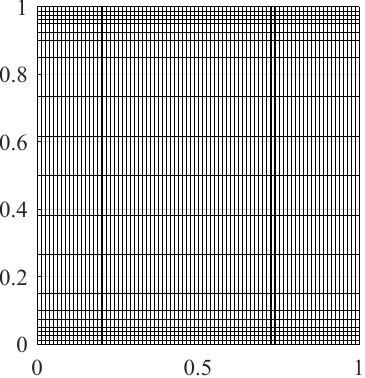}}
    
    \caption{Computational meshes for the subdomains of the convection-diffusion equation.}
    \label{fig:meshSize}
\end{figure}

The spatial domain is discretized using a structured grid consisting of 540 and 1,600 quadrilateral elements in $\Omega_1$ and $\hOmega_2$, respectively. A non-uniform mesh refinement is performed a priori in the vicinity of the walls, as detailed in Fig.~\ref{fig:meshSize}.
The parametric domains $\mathcal{I}_1$ and $\mathcal{I}_2$ are discretized using $10^{5}$ and $3.5 \times 10^{3}$ collocation points uniformly distributed, leading to mesh sizes in the parametric directions equal to $h_{\mu_1} = 10^{-1}$ and $h_{\mu_2} = 10^{-3}$.
Note that for all the values of the parameters under analysis, the P\'{e}clet number is larger than 1 and problem~\eqref{eq:paramConvDiff} is convection-dominated. Hence, the PGD-based surrogate model is constructed using a continuous Galerkin formulation with streamline upwind Petrov-Galerkin (SUPG) stabilization and $\mathbb{Q}_1$ Lagrange basis functions for the spatial problem, see~\cite{Huerta-GCCDH-13}.

The aforementioned spatial discretization yields interfaces with 19 DOFs.
The local surrogate model with dimensionality reduction is thus compared to the strategy proposed in~\cite{Discacciati:2024:CMAME}.
Table~\ref{tab:convDiffOffline} reports the details of the offline phase of the construction of the PGD-based surrogate model.
As previously observed, for each subdomain, the strategy in Sect.~\ref{sc:offlineOnline} requires the computation of $\nDP=1$ data problem and $\nIP=19$ interface problems. The latter feature $\dIP=3$ (2 spatial dimensions and 1 parameter) in $\Omega_1$ and $\dIP=4$ (2 spatial dimensions and 2 parameters) in $\hOmega_2$.
The approach with active interface parameters of Sect.~\ref{sc:pgdAIP} reduces the number of problems to be solved (7 in $\Omega_1$ and 10 in $\hOmega_2$), but increases their dimensionality to 6, with significantly worse performance of the offline phase.
It is worth noticing that, in order to balance the computational load of the computation of the surrogate model in the two subdomains, different numbers of active parameters are considered in $\Omega_1$ and $\hOmega_2$. In particular, to attain the same dimensionality in the two subdomains, $\nAIP=3$ is selected in $\Omega_1$ and $\nAIP=2$ in $\hOmega_2$.
\begin{table}[!htb]
\centering
\begin{tabular}{| l | c || c | c | c | c | c || c | c || c |}
\hline
Surrogate & \multirow{2}{*}{$\nAIP$} & \multirow{2}{*}{$\nDP$} & \multicolumn{2}{c|}{$\nIP$} & \multicolumn{2}{c||}{$\dIP$} & \multicolumn{2}{c||}{Modes$^\dagger$} &\multirow{2}{*}{$\Toff$ (s)} \\
\cline{4-9}
\rule{0pt}{11pt} strategy & & & $\Omega_1$ & $\hOmega_2$ & $\Omega_1$ & $\hOmega_2$ & $\Omega_1$ & $\hOmega_2$ & \\
\hline
Sect.~\ref{sc:offlineOnline} & - & 1 & 19 & 19 & 3 & 4 & 51 (75) & 157 (412) & 61.97 \\
\hline
\multirow{2}{*}{Sect.~\ref{sc:pgdAIP}} & 2 & 1 & - & 10 & - & 6 & - & 165 (596) & \multirow{2}{*}{6797.61} \\
& 3 & 1 & 7 & - & 6 & - & 38 (55) & - &  \\
\hline
\end{tabular}

\caption{Offline phase for the local surrogate models using two different approaches to handle the interface parameters. $^\dagger$The number in brackets denotes the number of modes before compression.}
\label{tab:convDiffOffline}
\end{table}

The overall number of modes computed after PGD compression using the model in Sect.~\ref{sc:offlineOnline} is 208, whereas the active interface parameters approach requires 203.
Nonetheless, before compression, the two methods determine 487 and 651 modes, respectively.
The lower number of computations and, most importantly, the reduced dimensionality of the interface problems ($\dIP$ being 3 or 4 instead of 6) are responsible for a significant improvement of performance with respect to the algorithm introduced in~\cite{Discacciati:2024:CMAME}, with a speed-up of approximately 110 times in the CPU time $\Toff$.

\begin{figure}[bht]
    \centering
    \subfigure[Surrogate with reduced dimensionality]{\includegraphics[width=0.7\columnwidth]{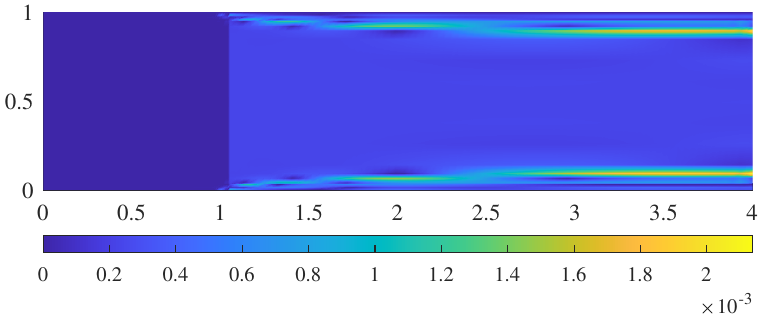}}
    
    \subfigure[Surrogate with $\nAIP=3$ in $\Omega_1$ and $\nAIP=2$ in $\hOmega_2$]{\includegraphics[width=0.7\columnwidth]{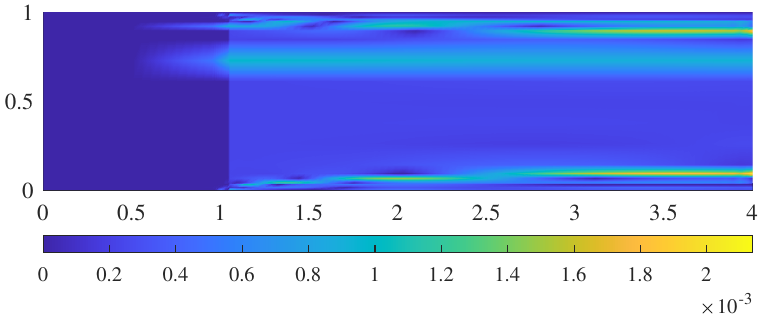}}
    
    \caption{Map of the scaled nodal error $| u_{\pgd}(\bmu) - u^h_{\Omega}(\bmu)|/\max_{\Omega} |u^h_{\Omega}(\bmu)|$ for $\bmu = (1.25 \times 10^4, 3)$ using the two different approaches to handle the interface parameters.}
    \label{fig:convDiffError}
\end{figure}

Finally, the accuracy of the two approaches is studied comparing the outcome of the online phase for $\bmu = (1.25 \times 10^4, 3)$.
Let $u^h_{\Omega}(\bmu)$ denote the corresponding finite element approximation of  problem~\eqref{eq:paramConvDiff}. Figure~\ref{fig:convDiffError} displays the map of the error $| u_{\pgd}(\bmu) - u^h_{\Omega}(\bmu)|/\max_{\Omega} |u^h_{\Omega}(\bmu)|$ of the PGD-based local surrogate model with reduced dimensionality and active interface parameters.
Whilst the overall magnitude of the error is comparable, the reduced dimensionality strategy removes numerical artifacts (most likely caused by the interpolation of the interface parameters) clearly visible near the interface and in the second domain when the method introduced in~\cite{Discacciati:2024:CMAME} is employed.
The GMRES iterations required to achieve convergence lower from 8 in the case of active interface parameters to 6 in the reduced dimensionality approach of Sect.~\ref{sc:offlineOnline}.
The corresponding online computing time $\Ton$ is also improved, performing the online coupling in $3 \times 10^{-2}~\text{s}$, roughly 5 times faster than the online phase timing of $1.45 \times 10^{-1}~\text{s}$ reported in~\cite{Discacciati:2024:CMAME}.

\subsection{Computational comparison of the online phase}
\label{sc:compareOnline}

In this section, the local surrogate models with reduced dimensionality and with active interface parameters are compared in terms of their performance during the online coupling procedure.
To this end, we consider a multi-domain benchmark test, featuring a two-dimensional thermal problem with discontinuous conductivities and 9 parameters, introduced in~\cite{Eftang:2013:IJNME}.

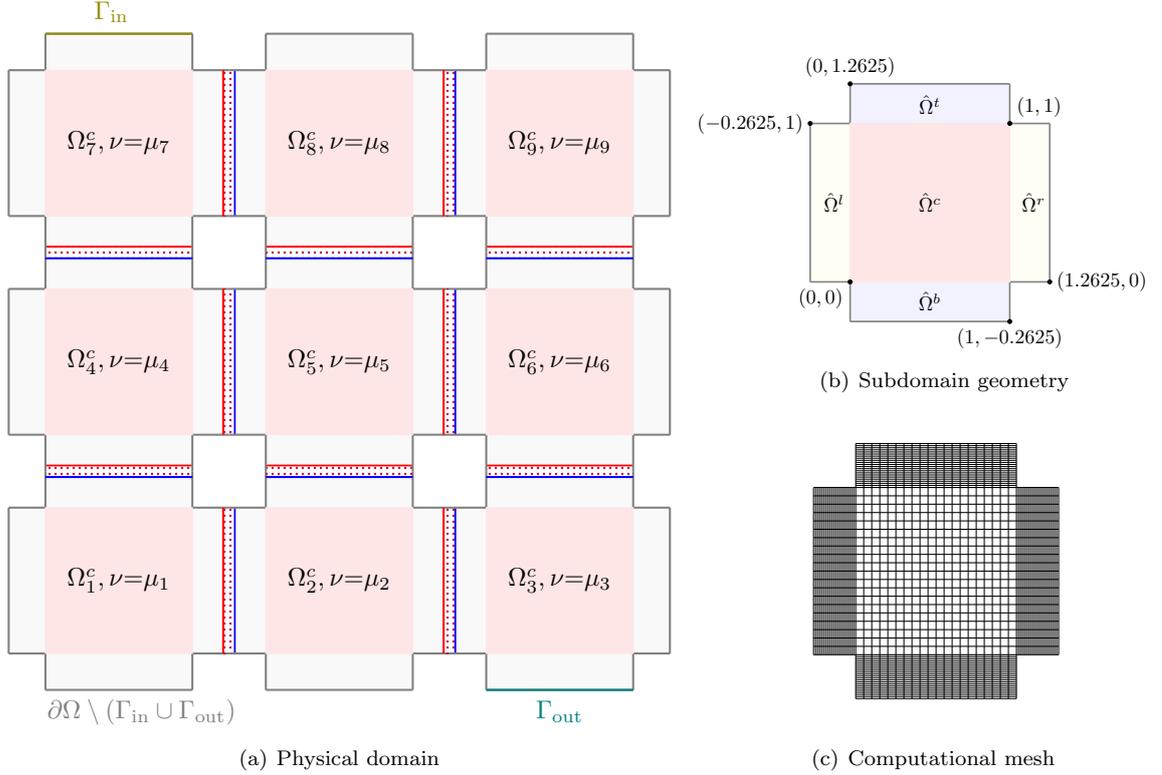
\begin{figure}[!htb]
	\centering
    \subfigure[Physical domain]{
		\resizebox{0.6\columnwidth}{!}{
			\begin{tikzpicture}
				\def\s{2}
				\begin{scope}
					\pateraWing{1}{0}{1}{0}
					\pateraBlock{1}{1}{1}{2}{2}{1}
				\end{scope}
				
				\begin{scope}[shift={(1.5*\s, 0)}]
					\pateraWing{1}{0}{0}{0}
					\pateraBlock{2}{1}{2}{2}{2}{1}
				\end{scope}
				
				\begin{scope}[shift={(3*\s, 0)}]
					\pateraWing{1}{0}{0}{1}
					\pateraBlock{3}{1}{2}{1}{2}{1}
					\draw[teal, thick] (0, -0.25*\s) -- (1*\s, -0.25*\s);
				\end{scope}
				
				\begin{scope}[shift={(0, 1.5*\s)}]
					\pateraWing{0}{0}{1}{0}
					\pateraBlock{4}{2}{1}{2}{2}{1}
				\end{scope}
				
				\begin{scope}[shift={(1.5*\s, 1.5*\s)}]
					\pateraWing{0}{0}{0}{0}
					\pateraBlock{5}{2}{2}{2}{2}{1}
				\end{scope}
				
				\begin{scope}[shift={(3*\s, 1.5*\s)}]
					\pateraWing{0}{0}{0}{1}
					\pateraBlock{6}{2}{2}{1}{2}{1}
				\end{scope}
				
				\begin{scope}[shift={(0, 3*\s)}]
					\pateraWing{0}{1}{1}{0}
					\pateraBlock{7}{2}{1}{2}{1}{1}
					\draw[olive, thick] (0, 1.25*\s) -- (1*\s, 1.25*\s);
				\end{scope}
				
				\begin{scope}[shift={(1.5*\s, 3*\s)}]
					\pateraWing{0}{1}{0}{0}
					\pateraBlock{8}{2}{2}{2}{1}{1}
				\end{scope}
				
				\begin{scope}[shift={(3*\s, 3*\s)}]
					\pateraWing{0}{1}{0}{1}
					\pateraBlock{9}{2}{2}{1}{1}{1}
				\end{scope}
				
				\node[olive] at (0.45*\s, 4.4*\s) {$\Gamma_{\text{in}}$};
				
				\node[teal] at (3.5*\s, -0.4*\s) {$\Gamma_{\text{out}}$};

                \node[gray, anchor=west] at (-0.05*\s, -0.4*\s) {$\partial\Omega \setminus (\Gamma_{\text{in}} \cup \Gamma_{\text{out}})$};
				
			\end{tikzpicture}
            \label{fig:pateraDom}
		}
	}
    \hspace{5pt}
    \begin{minipage}[b]{0.35\columnwidth}
    \subfigure[Subdomain geometry]{
        \hspace{-20pt}
        \resizebox{1.2\columnwidth}{!}{
        \begin{tikzpicture}
            \def\s{2.8}
            
            \pateraWing{2}{2}{2}{2}
            \pateraBlock{}{1}{1}{1}{1}{3}
    
            \filldraw[black] (0, 0) circle (1pt) node[anchor=north east]{$(0, 0)$};
            \filldraw[black] (1*\s, 1*\s) circle (1pt) node[anchor=south west]{$(1, 1)$};
    
            \filldraw[black] (1.25*\s, 0) circle (1pt) node[anchor=west]{$(1.2625, 0)$};
            \filldraw[black] (1*\s, -0.25*\s) circle (1pt) node[anchor=north]{$(1, -0.2625)$};
            \filldraw[black] (-0.25*\s, 1*\s) circle (1pt) node[anchor=east]{$(-0.2625, 1)$};
            \filldraw[black] (0, 1.25*\s) circle (1pt) node[anchor=south]{$(0, 1.2625)$};
        \end{tikzpicture}
        }
        \label{fig:subdomGeo}
    }

    \hspace{25pt}
    \subfigure[Computational mesh]{
      \includegraphics[width=0.65\columnwidth]{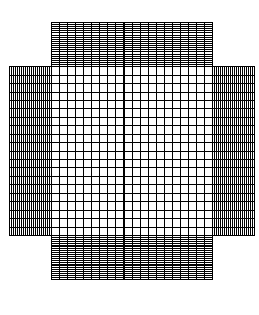}
      \label{fig:subdomMesh}
    }
  \end{minipage}
    
	\caption{Physical and computational domains for the multi-domain thermal problem with parametrized, discontinuous conductivities.}
	\label{fig:testPatera}
\end{figure}

Let $\Omega$ be the domain reported in Fig.~\ref{fig:pateraDom}.
Following the previously introduced color notation, blue and red lines identify the internal interfaces $\Gamma_1$ and $\Gamma_2$ for each pair of overlapping subdomains and the thin dotted purple regions denote the overlaps.
The parametric problem~\eqref{eq:globalProblemA} is characterized by $f(\bmu)=0$ and is equipped with following set of boundary conditions
\begin{equation}\label{eq:BCpatera}
\begin{aligned}
    u(\bmu) &= 0  && \text{on $\Gamma_{\text{out}}$,} \\
    \nu(\bmu) \nabla u(\bmu) \cdot \bn &= 1 && \text{on $\Gamma_{\text{in}}$,} \\
    \nu(\bmu) \nabla u(\bmu) \cdot \bn &= 0 && \text{on $\partial\Omega \setminus (\Gamma_{\text{in}} \cup \Gamma_{\text{out}})$,}
\end{aligned}    
\end{equation}
where $\bn$ denotes the outward unit normal to the corresponding boundary.
The thermal conductivity is defined by means of $N_p=9$ scalar parameters $\mu_k \in \mathcal{I}_k$, with $\mathcal{I}_k = [5\times 10^{-2},10]$, as
\begin{equation}\label{eq:conductivity}
\nu(\bmu) = \begin{cases}
\mu_i & \text{in $\Omega_i^c, \ i =1,\ldots,9$}, \\
1 & \text{otherwise} ,
\end{cases}
\end{equation}
such that a unique constant value equal to 1 is introduced in the overlapping regions, whereas the conductivity inside the central region $\Omega_i^c$ of each subdomain is parametrized (see Fig.~\ref{fig:testPatera}).

The physical domain $\Omega$ is subdivided in 9 overlapping subdomains $\Omega_i, \ i=1,\ldots,9$, with overlaps located in the bottom, top, and lateral wings, of dimension $2.5 \times 10^{-2}$.
The geometry of each subdomain is detailed in Fig.~\ref{fig:subdomGeo} and the corresponding computational grid of 2,080 quadrilaterals is displayed in Fig.~\ref{fig:subdomMesh}.
Let $\hOmega^c$ denote the central region (pink), $\hOmega^b$ and $\hOmega^t$ the bottom and top wings (blue), and $\hOmega^l$ and $\hOmega^r$ the left and right wings (yellow).
A non-uniform grid with varying horizontal, $h_x$, and vertical, $h_y$, mesh sizes is defined, namely,
\begin{equation}
h_x =
\begin{cases}
5 \times 10^{-2} & \text{in $\hOmega^c \cup \hOmega^b \cup \hOmega^t$,}\\
1.25 \times 10^{-2}  & \text{in $\hOmega^l \cup \hOmega^r$,}
\end{cases}
\qquad
h_y =
\begin{cases}
5 \times 10^{-2} & \text{in $\hOmega^c \cup \hOmega^l \cup \hOmega^r$,}\\
1.25 \times 10^{-2}  & \text{in $\hOmega^b \cup \hOmega^t$.}
\end{cases}
\end{equation}

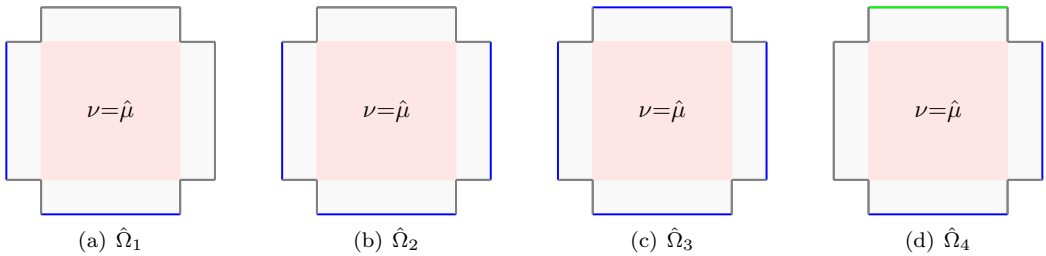
\begin{figure}[h!]
    \centering
    \subfigure[$\hOmega_1$]{
        \resizebox{0.2\textwidth}{!}{
        \begin{tikzpicture}
            \def\s{2}
            \pateraWing{1}{1}{1}{1}
            \pateraBlock{1}{0}{0}{1}{1}{2}

            \draw[blue, thick] (-0.25*\s, 0) -- (-0.25*\s, 1*\s);

            \draw[blue, thick] (0, -0.25*\s) -- (1*\s, -0.25*\s);

            \draw[gray, thick] (0, 1.25*\s) -- (1*\s, 1.25*\s);
            \draw[gray, thick] (1.25*\s, 0) -- (1.25*\s, 1*\s);
            \draw[gray, thick] (-0.25*\s, 1*\s) -- (0, 1*\s);
            \draw[gray, thick] (1*\s, 1*\s) -- (1.25*\s, 1*\s);
            \draw[gray, thick] (-0.25*\s, 0) -- (0, 0);
            \draw[gray, thick] (1*\s, 0) -- (1.25*\s, 0);
            \draw[gray, thick] (0, -0.25*\s) -- (0, 0);
            \draw[gray, thick] (0, 1*\s) -- (0, 1.25*\s);
            \draw[gray, thick] (1*\s, -0.25*\s) -- (1*\s, 0);
            \draw[gray, thick] (1*\s, 1*\s) -- (1*\s, 1.25*\s);
        \end{tikzpicture}
        }
    }
    \hspace{5pt}
    \subfigure[$\hOmega_2$]{
        \resizebox{0.2\textwidth}{!}{
        \begin{tikzpicture}
            \def\s{2}
            \pateraWing{1}{1}{1}{1}
            \pateraBlock{2}{0}{0}{0}{1}{2}

            \draw[blue, thick] (-0.25*\s, 0) -- (-0.25*\s, 1*\s);

            \draw[blue, thick] (0, -0.25*\s) -- (1*\s, -0.25*\s);

            \draw[blue, thick] (1.25*\s, 0) -- (1.25*\s, 1*\s);

            \draw[gray, thick] (0, 1.25*\s) -- (1*\s, 1.25*\s);
            \draw[gray, thick] (-0.25*\s, 1*\s) -- (0, 1*\s);
            \draw[gray, thick] (1*\s, 1*\s) -- (1.25*\s, 1*\s);
            \draw[gray, thick] (-0.25*\s, 0) -- (0, 0);
            \draw[gray, thick] (1*\s, 0) -- (1.25*\s, 0);
            \draw[gray, thick] (0, -0.25*\s) -- (0, 0);
            \draw[gray, thick] (0, 1*\s) -- (0, 1.25*\s);
            \draw[gray, thick] (1*\s, -0.25*\s) -- (1*\s, 0);
            \draw[gray, thick] (1*\s, 1*\s) -- (1*\s, 1.25*\s);
        \end{tikzpicture}
        }
    }
    \hspace{5pt}
    \subfigure[$\hOmega_3$]{
        \resizebox{0.2\textwidth}{!}{
        \begin{tikzpicture}
            \def\s{2}
            \pateraWing{1}{1}{1}{1}
            \pateraBlock{3}{0}{0}{0}{0}{2}

            \draw[blue, thick] (-0.25*\s, 0) -- (-0.25*\s, 1*\s);

            \draw[blue, thick] (0, -0.25*\s) -- (1*\s, -0.25*\s);

            \draw[blue, thick] (1.25*\s, 0) -- (1.25*\s, 1*\s);

            \draw[blue, thick] (0, 1.25*\s) -- (1*\s, 1.25*\s);

            %
            \draw[gray, thick] (-0.25*\s, 1*\s) -- (0, 1*\s);
            \draw[gray, thick] (1*\s, 1*\s) -- (1.25*\s, 1*\s);
            \draw[gray, thick] (-0.25*\s, 0) -- (0, 0);
            \draw[gray, thick] (1*\s, 0) -- (1.25*\s, 0);
            \draw[gray, thick] (0, -0.25*\s) -- (0, 0);
            \draw[gray, thick] (0, 1*\s) -- (0, 1.25*\s);
            \draw[gray, thick] (1*\s, -0.25*\s) -- (1*\s, 0);
            \draw[gray, thick] (1*\s, 1*\s) -- (1*\s, 1.25*\s);
        \end{tikzpicture}
       } 
    }
    \hspace{5pt}
    \subfigure[$\hOmega_4$]{
        \resizebox{0.2\textwidth}{!}{
        \begin{tikzpicture}
            \def\s{2}
            \pateraWing{1}{1}{1}{1}
            \pateraBlock{4}{0}{1}{0}{1}{2}

            \draw[blue, thick] (0, -0.25*\s) -- (1*\s, -0.25*\s);

            \draw[blue, thick] (1.25*\s, 0) -- (1.25*\s, 1*\s);

            \draw[green, thick] (0, 1.25*\s) -- (1*\s, 1.25*\s);

            \draw[gray, thick] (-0.25*\s, 0) -- (-0.25*\s, 1*\s);
            \draw[gray, thick] (-0.25*\s, 1*\s) -- (0, 1*\s);
            \draw[gray, thick] (1*\s, 1*\s) -- (1.25*\s, 1*\s);
            \draw[gray, thick] (-0.25*\s, 0) -- (0, 0);
            \draw[gray, thick] (1*\s, 0) -- (1.25*\s, 0);
            \draw[gray, thick] (0, -0.25*\s) -- (0, 0);
            \draw[gray, thick] (0, 1*\s) -- (0, 1.25*\s);
            \draw[gray, thick] (1*\s, -0.25*\s) -- (1*\s, 0);
            \draw[gray, thick] (1*\s, 1*\s) -- (1*\s, 1.25*\s);
        \end{tikzpicture}
        }
    }

    \caption{Reference subdomains $\hOmega_j, \ j=1,\ldots,4$. Boundary condition type: Dirichlet (blue), homogeneous Neumann (grey), non-homogeneous Neumann (green).}
    \label{fig:pateraRef}
\end{figure}

Following~\cite{Discacciati:2024:CMAME}, the computational domain can be constructed by composing four reference subdomains $\hat{\Omega}_j, \ j=1,\ldots,4$ (see Fig.~\ref{fig:pateraRef}) with suitable rigid rotations and/or translations (see Table~\ref{tab:pateraTransform}). A new parameter $\hmu$ is introduced and, for each reference subdomain, the conductivity is thus defined in terms of a unique parameter as
\begin{equation}\label{eq:conductivityRef}
\nu(\hmu) = \begin{cases}
\hmu & \text{in $\hOmega^c$}, \\
1 & \text{in $\hOmega^b \cup \hOmega^t \cup \hOmega^l \cup \hOmega^r$} ,
\end{cases}
\end{equation}
with $\hmu \in [5\times 10^{-2},10]$.
The spatial problem is approximated using a continuous Galerkin finite element method with $\mathbb{Q}_1$ Lagrange basis functions, whereas the parametric interval is discretized with pointwise collocation on a uniform grid with $9.95 \times 10^3$ equally-spaced nodes and $h_{\hmu} = 10^{-3}$.
The resulting number of DOFs of the finite element approximation and the corresponding number of unknowns on the interfaces is reported in Table~\ref{tab:pateraDOFs} for each subdomain.
\begin{table}[!htb]
\begin{center}
\resizebox{\textwidth}{!}{%
\begin{tabular}{| l || c | c | c | c | c | c | c | c | c|}
\hline
Phys. subdomain & $\Omega_1$  & $\Omega_2$  & $\Omega_3$  & $\Omega_4$  & $\Omega_5$  & $\Omega_6$  & $\Omega_7$  & $\Omega_8$  & $\Omega_9$ \\
Conductivity $\hmu$ & 0.1 & 0.2 & 0.4 & 0.8 & 1.6 & 3.2 & 6.4 & 0.1 & 0.2 \\
\hline
\rule{0pt}{11pt}Ref. subdomain & $\hOmega_1$ & $\hOmega_2$ & $\hOmega_2$ & $\hOmega_2$ & $\hOmega_3$ & $\hOmega_2$ & $\hOmega_4$ & $\hOmega_2$ & $\hOmega_1$ \\
\rule{0pt}{18pt}Translation & $\begin{pmatrix}0 \\ 0\end{pmatrix}$ & $\begin{pmatrix} 1.5 \\ 0\end{pmatrix}$ & $\begin{pmatrix}3 \\ 0\end{pmatrix}$ & $\begin{pmatrix}0 \\ 1.5 \end{pmatrix}$ & $\begin{pmatrix} 1.5 \\ 1.5 \end{pmatrix}$ & $\begin{pmatrix} 3 \\ 1.5\end{pmatrix}$ & $\begin{pmatrix} 0 \\ 3\end{pmatrix}$ & $\begin{pmatrix} 1.5 \\ 3 \end{pmatrix}$ & $\begin{pmatrix} 3 \\ 3 \end{pmatrix}$ \\
\rule{0pt}{18pt}Rotation    & $\pi$ & $\pi$ & $\displaystyle\frac{3\pi}{2}$ & $\displaystyle\frac{\pi}{2}$ & 0 & $\displaystyle\frac{3\pi}{2}$ & 0 & 0 & 0\\[0.5em]
\hline
\end{tabular}}
\end{center}
\caption{Transformations of the reference subdomains $\hOmega_j, \ j=1,\ldots,4$ into the physical subdomains $\Omega_i, \ i=1,\ldots,9$.}
\label{tab:pateraTransform}
\end{table}
\begin{table}[!htb]
\begin{center}
\begin{tabular}{| c | c || c | c |}
\hline
Reference & FEM & Physical & Interface \\
subdomain & DOFs & parameters & parameters \\
\hline
\rule{0pt}{11pt}$\hOmega_1$ & 2,163 & 1 & 42\phantom{$^\dagger$} \\
$\hOmega_2$ & 2,142 & 1 & 63$^\dagger$ \\
$\hOmega_3$ & 2,121 & 1 & 84\phantom{$^\dagger$} \\
$\hOmega_4$ & 2,163 & 1 & 42\phantom{$^\dagger$} \\
\hline
\end{tabular}
\end{center}

\caption{Dimensions of the local subproblems. $^\dagger$In subdomain $\Omega_3$, a homogeneous Dirichlet condition is enforced on the boundary $\Gamma_{\text{out}}$. Hence, of the 63 interface parameters of $\hOmega_2$, 21 nodes are fixed a priori and only 42 remain to be determined by the coupling procedure}.
\label{tab:pateraDOFs}
\end{table}

The details of the offline phase are reported in Table~\ref{tab:pateraOffline}, where the number of subproblems and their dimensionality for the local surrogate model with reduced dimensionality (Sect.~\ref{sc:offlineOnline}) are compared to those of the active interface parameters approach introduced in~\cite{Discacciati:2024:CMAME} with $\nAIP=3$.
As observed in Sect.~\ref{sc:sensitivityTesting} and~\ref{sc:compareOffline}, the method in Sect.~\ref{sc:offlineOnline} clearly outperforms the local surrogate model based on active interface parameters: the lower dimension $\dIP$ of the local subproblems allows the PGD to converge significantly faster, computing fewer modes and achieving a speed-up of approximately 24 times.
\begin{table}[!htb]
\centering
\begin{tabular}{| l | c || c | c | c | c | c | c |}
\hline
Surrogate & \multirow{2}{*}{$\nAIP$} & Reference & \multirow{2}{*}{$\nDP$} & \multirow{2}{*}{$\nIP$} & \multirow{2}{*}{$\dIP$} & \multirow{2}{*}{Modes$^\dagger$} & \multirow{2}{*}{$\Toff$ (s)} \\
strategy & & subdomain & & & & & \\
\hline
\rule{0pt}{11pt}\multirow{4}{*}{Sect.~\ref{sc:offlineOnline}} & \multirow{4}{*}{-} & $\hOmega_1$ & 0 & 42 & 3 & 190 (302) & \multirow{4}{*}{36.18} \\
& & $\hOmega_2$ & 0 & 63 & 3 & 272 (456) & \\
& & $\hOmega_3$ & 0 & 84 & 3 & 336 (604) & \\
& & $\hOmega_4$ & 1 & 42 & 3 & 195 (311) & \\
\hline
\rule{0pt}{11pt}\multirow{4}{*}{Sect.~\ref{sc:pgdAIP}} & \multirow{4}{*}{3} & $\hOmega_1$ & 0 & 14 & 6 & 236 (370) & \multirow{4}{*}{872.38} \\
& & $\hOmega_2$ & 0 & 21 & 6 & 386 (672) & \\
& & $\hOmega_3$ & 0 & 28 & 6 & 509 (930) & \\
& & $\hOmega_4$ & 1 & 14 & 6 & 244 (419) & \\
\hline
\end{tabular}

\caption{Offline phase for the local surrogate models using two different approaches to handle the interface parameters. $^\dagger$The number in brackets denotes the number of modes before compression.}
\label{tab:pateraOffline}
\end{table}

\begin{figure}[!htb]
    \centering
    \subfigure[Surrogate with reduced dimensionality]{\includegraphics{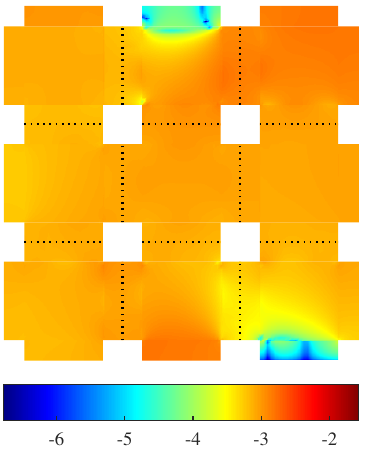}}
    \hspace{10pt}
    \subfigure[Surrogate with $\nAIP=3$]{\includegraphics{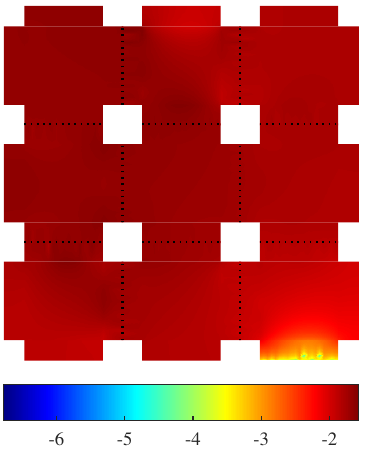}}
    
    \caption{Map of the scaled nodal error $\log_{10}(| u_{\pgd}(\hmu) - u^h(\hmu)|/\max_{\Omega} |u^h(\hmu)|)$ for $\hmu$ reported in Table~\ref{tab:pateraTransform} using two different approaches to handle the interface parameters. The dashed lines indicate the location of the overlapping regions.}
    \label{fig:pateraComparison}
\end{figure}

To compare the performance of the online phase, we consider the set of parameters described in Table~\ref{tab:pateraTransform}.
Let $u^h(\hmu)$ denote the solution obtained with a high-fidelity overlapping DD-FEM approximation. Figure~\ref{fig:pateraComparison} displays the map of the error $\log_{10}(| u_{\pgd}(\hmu) - u^h(\hmu)|/\max_{\Omega} |u^h(\hmu)|)$ of the PGD-based local surrogate model with reduced dimensionality and active interface parameters.
It is straightforward to observe that the surrogate model with reduced dimensionality provides one order of magnitude extra accuracy with respect to the strategy based on clustering the interface nodes. Moreover, the results confirm that the coupling algorithm is extremely accurate even using a surrogate model, with no significant error being introduced at the subdomain interfaces.
Besides the qualitative improvement observed in Fig.~\ref{fig:pateraComparison}, Table~\ref{tab:pateraOnlineRes} reports the relative error measured in $\ell^{\infty}(\Omega)$ norm for a fixed value of the parameter $\hmu$, namely,
\begin{equation}\label{eq:errInfDef}
E^{\pgd}_{\infty} = \frac{\| u_{\pgd}(\hmu) - u^h(\hmu) \|_{\ell^{\infty}(\Omega)}}{\| u^h(\hmu) \|_{\ell^{\infty}(\Omega)}} .
\end{equation}
The local surrogate model with reduced dimensionality outperforms the strategy based on active interface parameters, lowering the relative error $E^{\pgd}_{\infty}$ by one order of magnitude and achieving an accuracy of $10^{-3}$.

\begin{table}[!htb]
\centering
\begin{tabular}{| l | c || c | c | c |}
\hline
Surrogate & \multirow{2}{*}{$\nAIP$} & \multirow{2}{*}{$\nGMRES$} & \multirow{2}{*}{$E^{\pgd}_{\infty}$} & \multirow{2}{*}{$\Ton$ (s)} \\
strategy & & & & \\
\hline
\rule{0pt}{11pt}Sect.~\ref{sc:offlineOnline} & - & 93 & $1.8 \times 10^{-3}$ & 0.47 \\
\hline
\rule{0pt}{11pt}Sect.~\ref{sc:pgdAIP} & 3 & 302 & $2.5 \times 10^{-2}$ & 157.23 \\
\hline\hline
\multicolumn{2}{|l||}{\rule{0pt}{11pt}Reference DD-FEM} & 95 & - & 414.80 \\
\hline
\end{tabular}

\caption{Online phase for the local surrogate models using two different approaches to handle the interface parameters.}
\label{tab:pateraOnlineRes}
\end{table}

The number of GMRES iterations required to attain convergence is reduced from 302 to 93. The fewer iterations are likely due to the avoidance of the interpolation procedure performed by the clustering approach.
This is particularly significant since the reference full-order DD-FEM method requires 95 GMRES iterations to converge, thus the local surrogate model with reduced dimensionality is capable of providing performance comparable to the high-fidelity solver.
In addition, this can be done in real time, with an online execution time $\Ton$ of less than half a second, with a speed-up of
334 times with respect to the active interface parameters case and of 882 with respect to DD-FEM.
The CPU times reported in the table include the precomputations to setup the iterative solver and solve the interface system of dimension 504.

Similar results, not reported here for brevity, are also obtained for another test case studied in~\cite{Discacciati:2024:CMAME}, where the variations of the conductivity parameter are less pronounced.
Also in this case, the local surrogate model with reduced dimensionality requires a number of GMRES iterations comparable to the high-fidelity DD-FEM, while achieving errors of the order of $10^{-4}$ and significantly reducing the CPU time of the previously proposed strategy based on active interface parameters.

\section{Concluding remarks}
\label{sc:Conclusion}
In this work, we presented a novel approach to construct physics-based local surrogate models based on the proper generalized decomposition using overlapping domain decomposition.
The method starts from the framework introduced in~\cite{Discacciati:2024:CMAME} and significantly improves its performance, both in terms of accuracy and efficiency.

DD-PGD exploits the linearity of the parametric PDE under analysis to define low-dimensional local subproblems with arbitrary Dirichlet boundary conditions at the interface.
Whilst the approach in~\cite{Discacciati:2024:CMAME} clusters the interface nodes, the surrogate model with reduced dimensionality proposed in the present work introduces a novel definition of the trace variable to devise local subproblems with unitary boundary conditions at each node of the interface.
This yields a set of subproblems with the same number of spatial and parametric dimensions as the original parametric equation, circumventing the challenge of accurately and efficiently describing the trace of the solution, commonly experienced by ROM-based domain decomposition algorithms.

The method inherits the advantages of the original DD-PGD framework, including: (i) non-intrusiveness with respect to the underlying full-order solver, since in the offline phase it only requires the imposition of Dirichlet boundary conditions at any node on the interface using the traces of the finite element basis functions; (ii) seamless implementation, since it does not require the introduction of Lagrange multipliers or additional variables to impose the continuity of the solution in the overlap; (iii) real-time coupling, since in the online phase the method only solves a linear system of equations for the nodal values of the solution at the interface with no additional parametric problem to be solved.

The work performs a detailed computational study of the performance of the local surrogate model with reduced dimensionality, comparing it with the algorithm based on active interface parameters presented in~\cite{Discacciati:2024:CMAME}.
Three benchmarks are presented: (i) a sensitivity study of the effect of the number of parametric dimensions on the solution of the local subproblems using a two-domain parametric Poisson equation; (ii) a computational assessment of the offline phase using a convection-dominated convection-diffusion equation with parametrized geometry; (iii) a computational evaluation of the online phase via a multi-domain parametric thermal problem.
The local surrogate model with reduced dimensionality outperforms the DD-PGD algorithm in~\cite{Discacciati:2024:CMAME} by significantly accelerating both the offline and the online phase.
In the offline phase, reducing the number of dimensions of the local subproblems by 2 or 3 allows reducing the number of computed modes by approximately $25\%$, while achieving speed-ups up to 110 times.
In the online phase, the proposed approach achieves convergence with a number of GMRES iterations comparable to the high-fidelity DD-FEM method, reducing the number of iterations performed by the active interface parameter scheme by $69\%$. The corresponding CPU time is reduced from $157$~s to less than half a second, achieving real-time evaluation capabilities with a speed-up of 334 times.

\paragraph{Acknowledgements} The authors acknowledge funding as follows. MD: EPSRC grant EP/V027603/1. BJE: EPSRC Doctoral Training Partnership grant EP/W523987/1. MG: Spanish Ministry of Science, Innovation and Universities and Spanish State Research Agency \\ MICIU/AEI/10.13039/501100011033 (Grant No. PID2023-149979OB-I00); Generalitat de Catalunya
(Grant No. 2021-SGR-01049); MG is Fellow of the Serra H\'unter Programme of the Generalitat de Catalunya.


\bibliographystyle{elsarticle-num}
\bibliography{references_paper_ellipticUnitary}



	
\end{document}